\pdfoutput=1
\documentclass{shinyart}

\usepackage[utf8]{inputenc}

\usepackage{shinybib}

\usepackage{enumitem}
\setlist[enumerate]{label=(\roman*)}

\usepackage{autonum}

\newcommand{\eps}{\varepsilon}
\renewcommand{\phi}{\varphi}
\newcommand{\N}{\mathbb{N}}
\newcommand{\R}{\mathbb{R}}

\newcommand{\calI}{\mathcal{I}}

\newcommand{\calU}{\mathcal{U}}

\newcommand{\norm}[1]{\| #1 \|}

\DeclareMathOperator{\sign}{\mathrm{sign}}

\newcommand{\prox}{\mathrm{prox}}
\newcommand{\proj}{\mathrm{proj}}
\newcommand{\swerror}{\sigma_{\mathrm{sw}}}
 
\newcommand{\omobs}{\omega_\mathrm{obs}}

\usepackage{booktabs}
\usepackage{tikz,pgfplots}
\pgfplotsset{compat=newest}
\pgfplotsset{plot coordinates/math parser=false}

\title{Nonconvex penalization of switching control of partial differential equations}

\author{Christian Clason\thanks{Faculty of Mathematics, University Duisburg-Essen, 45117 Essen, Germany (\email{christian.clason@uni-due.de})}
    \and Karl Kunisch\thanks{Institute of Mathematics and Scientific Computing, University of Graz, Heinrichstrasse 36, 8010 Graz, Austria and Radon Institute, Austrian Academy of Sciences, Linz, Austria (\email{karl.kunisch@uni-graz.at})}
    \and Armin Rund\thanks{Institute of Mathematics and Scientific Computing, University of Graz, Heinrichstrasse 36, 8010 Graz, Austria (\email{armin.rund@uni-graz.at})}
}
\date{\today}

\hypersetup{
    pdftitle={Nonconvex penalization of switching control of partial differential equations},
    pdfauthor={Christian Clason, Karl Kunisch, Armin Rund},
    pdfkeywords={optimal control, switching control, partial differential equations, nonsmooth optimization, nonconvex analysis, semi-smooth Newton method}
}

\begin{document}

\maketitle

\begin{abstract}
    This paper is concerned with optimal control problems for parabolic partial differential equations with pointwise in time switching constraints on the control.
    A standard approach to treat constraints in nonlinear optimization is penalization, in particular using $L^1$-type norms. Applying this approach to the switching constraint leads to a nonsmooth and nonconvex infinite-dimensional minimization problem which is challenging both analytically and numerically.
    Adding $H^1$ regularization or restricting to a finite-dimensional control space allows showing existence of optimal controls. First-order necessary optimality conditions are then derived using tools of nonsmooth analysis. Their solution can be computed using a combination of Moreau--Yosida regularization and a semismooth Newton method. Numerical examples illustrate the properties of this approach.
\end{abstract}

\section{Introduction}

Switching control refers to time-dependent optimal control problems with a vector-valued control of which at most one component should be active at any point in time.
To partially set the stage, we consider for example optimal tracking control for a linear evolution equation
$y_t+Au = Bu$ on $\Omega_T:=(0,T]\times \Omega$ together with initial conditions $y(0)=y_0$ on $\Omega$, where $A$ is a linear second order elliptic operator defined on $\Omega\subset \R^n$ with homogeneous Neumann boundary conditions and the linear  control operator $B:L^2(0,T;\R^N)\to L^2(\Omega_T)$ is given by
\begin{equation}\label{eq:kk1}
    (Bu)(t,x)= \sum_{i=1}^N\chi_{\omega_i}(x) u_i(t),
\end{equation}
where $\chi_{\omega_i}$ are the characteristic functions of  given control domains $\omega_i\subset \Omega$ of positive measure.
Furthermore, let $\omobs\subset \Omega$ denote the observation domain and let $y^d\in L^2(0,T;L^2(\omobs))$ denote the target.
Consider now the standard optimal control problem
\begin{equation}\label{eq:problem_linquad}
    \left\{\begin{aligned}
            \min_{u \in L^2(0,T;\R^N)} &\frac12\norm{y-y^d}_{L^2(0,T;L^2(\omobs))}^2 + \frac\alpha2 \int_0^T|u(t)|_2^2\,dt,\\
            \text{s.\,t.}\quad & y_t+Ay = Bu, \quad y(0)=y_0,
    \end{aligned}\right.
\end{equation}
where $|v|_2^2 =\sum_{j=1}^N v_j^2$ denotes the squared $\ell^2$-norm on $\R^N$.
To promote the switching structure of optimal controls, we suggest adding the penalty term
\begin{equation}\label{eq:L1}
    \beta\int_0^T\sum_{\stackrel{i,j=1}{i< j}}^N |u_i(t) u_j(t)|\,dt
\end{equation}
with $\beta>0$ to the objective, which can be interpreted as an $L^1$-penalization of the switching constraint $u_i(t)u_j(t)=0$ for $i\neq j$ and $t\in[0,T]$.
The combination of control cost and switching penalty is convex if and only if $\beta\leq\alpha$.
The case $\beta=\alpha$ was investigated in \cite{switchingcontrol}; the aim of this work is to treat the case $\beta>\alpha$, which allows choosing the switching penalty parameter independently of the control cost parameter. As can be verified for a simple scalar example, there exist sets of data for which the minimizer of the convex problem is not switching, while the nonconvex problem does admit (possibly multiple) minimizers that are switching.

In the nonconvex case, the approach followed in \cite{switchingcontrol} is not applicable.
The main difficulty stems from the fact that the integrand $g:\R^2\to \R$, $(u_1,u_2)\mapsto |u_1u_2|$, is not convex, and hence the integral functional $G:L^2(0,T;\R^2)\to \R$, $u\mapsto \int_0^T |u_1(t)u_2(t)|\,dt$, is not weakly lower semicontinuous, which is an obstacle for proving existence. It is therefore necessary to enforce strong convergence of minimizing sequences, which is possible by either considering piecewise constant and hence finite-dimensional controls or by introducing an additional (small) $H^1(0,T;\R^2)$ penalty. Our analysis will cover both approaches.
Besides the question of existence of optimal controls, their numerical computation is also challenging due to the nonconvexity of the problem. Here we proceed as follows: Using the calculus of Clarke's generalized derivative \cite{Clarke:1990a,Clarke:2013}, we can derive first-order necessary optimality conditions. It then suffices to apply a Moreau--Yosida regularization only to the nonsmooth but convex term in the optimality conditions in order to apply a semismooth Newton method.

This is a natural continuation of our previous works \cite{CIK:2014,switchingcontrol} on convex relaxation of the switching constraint.
Let us briefly remark on further related literature. On switching control of ordinary and partial differential equations, there exists a large body of work; here we only mention \cite{Dolcetta:1984,Liberzon:2003,Shorten:2007} in the former context and \cite{Zuazua:2011,Lu:2013,Wang:2014,Gugat:2008,Hante:2013,Seidman} in the latter. A related topic is the control of switched systems, where we refer to, e.g., \cite{Hante:2009,Stojanovic:1989a,Stojanovic:1989b}.

This paper is organized as follows. \Cref{sec:existence} is concerned with existence of optimal controls and their convergence as $\beta\to \infty$ to a ``hard switching constrained'' problem. Optimality conditions are then derived in \cref{sec:optimality}, where the question of exact penalization is addressed as well. \Cref{sec:solution} discusses the numerical solution of the optimality conditions using a semismooth Newton method. Finally, \cref{sec:examples} presents numerical examples illustrating the properties of the nonconvex penalty approach.

\section{Existence}\label{sec:existence}

Here we describe the general framework that will be utilized and which will contain the example in the Introduction as a special case. Let $W$ denote a Hilbert space of measurable functions on the space-time cylinder $\Omega_T=(0,T]\times \Omega$, where $\Omega\subset \R^n$ is a bounded domain with Lipschitz continuous boundary. This space will serve as the state space of the solutions of the control system which appears as a constraint in \eqref{eq:problem_linquad}. It is assumed that $W \hookrightarrow L^2(0,T;L^2(\Omega))$, and that the embedding is continuous. Further let $\calU\subset L^2(0,T;\R^N)$ denote the Hilbert space of controls. We assume that there exists an affine control-to-state mapping $u\mapsto S(u)$. Here, we suppress the dependence of $S$ on $y_0$; for $y_0=0$, we denote the corresponding linear solution operator by $S_0$. Throughout it is assumed that $S$ satisfies
\begin{enumerate}[label=(\textsc{a}\arabic*),ref=\textsc{a}\arabic*]
    \item \label{A1}
        $S: L^2(0,T;\R^N) \to W$ is a continuous mapping satisfying
        \begin{equation}
            \norm{S(u)}_W \le C(\norm{u}_{L^2(0,T;\R^N)} + \norm{y_0}_{L^2(\Omega)})
        \end{equation}
        for a constant $C$ independent of $u$ and $y_0$.
\end{enumerate}
As mentioned in the Introduction, we need to restrict the set of feasible controls in order to obtain existence of an optimal control. We thus consider the following two cases for $\calU\subset L^2(0,T;\R^N)$:
\begin{enumerate}
    \item $\calU=H^1(0,T;\R^N)$;
    \item $\calU$ is finite-dimensional (e.g., consisting of piecewise constant controls).
\end{enumerate}

For the sake of presentation, we further restrict ourselves in the following to the case of two control components; the results remain valid for $N>2$ components (although it should be pointed out that, in contrast to the convex approach in \cite{switchingcontrol}, the number of terms in \eqref{eq:L1} grows as $\binom{N}{2}$).
We hence consider for $\beta>\alpha>0$ the problem
\begin{equation}\label{eq:problem}
    \min_{u \in \calU} \frac12\norm{Su-y^d}_{L^2(0,T;L^2(\omobs))}^2 + \frac\alpha2 \norm{u}_{L^2(0,T;\R^2)}^2+ \frac\eps2 \norm{u_t}_{L^2(0,T;\R^2)}^2 + \beta \int_0^T |u_1(t)u_2(t)|\,dt
\end{equation}
with $\omobs\subset\Omega$ and $y^d\in L^2(0,T;L^2(\omobs))$ as before.
If $\calU$ is finite-dimensional, it is understood that $\eps=0$; otherwise we require $\eps>0$.
Keeping $\eps\geq 0$ fixed, we will denote the cost functional in \eqref{eq:problem} by $J_\beta$.

\bigskip

Before we turn to address existence for \eqref{eq:problem}, we describe three typical cases of interest for which assumption \eqref{A1} is satisfied.  
Throughout the following, $A$ will denote a linear second-order uniformly elliptic operator with smooth coefficients.
\paragraph{Distributed control}
We return to the case considered in the Introduction, i.e.,
we consider the equation in \eqref{eq:problem_linquad} with $A$ together with homogenous Dirichlet, Neumann, or Robin boundary conditions and the control operator $B\in \mathcal{L}(L^2(0,T;\R^N), L^2(\Omega_T))$  as in \eqref{eq:kk1}. It is then well-known, see, e.g., \cite[Chap.~4]{Wloka}, that \eqref{A1} is satisfied with $W = W(0,T):=H^1(0,T;V^*)\cap L^2(0,T;V)$, where $V=H^1_0(\Omega)$ in the case of homogenous Dirichlet boundary conditions and $V=H^1(\Omega)$ for homogenous Neumann or Robin conditions.

\paragraph{Neumann boundary control}
Here we consider the case of Neumann boundary control. Thus the control system is given by
\begin{equation}
    \left\{
        \begin{aligned}
            y_t+Ay &= 0 \text{ in } Q_T, \\
            \frac{\partial y}{\partial n} &= Bu \text{ on } \Sigma_T, \\
            y(0)&=y_0 \text{ in } \Omega,
        \end{aligned}
    \right.
\end{equation}
where $\Sigma_T := (0,T]\times \partial \Omega$, and analogous to \eqref{eq:kk1} we now take $B$ to be of the form
\begin{equation}\label{eq:kk2}
    (Bu)(t,s)= \sum_{i=1}^2\chi_{\omega_i}(s) u_i(t),
\end{equation}
with $\chi_{\omega_i}$ the characteristic functions of  given control domains $\omega_i\subset \partial\Omega$ of positive measure relative to $\partial\Omega$. Again, \eqref{A1} is satisfied, this time with $W=W(0,T)$ and $V=H^1(\Omega)$. For a reference, see, e.g., \cite[Chap.~3.3]{Troeltzsch} and the references given there.

\paragraph{Dirichlet boundary control}
Finally we consider the case of Dirichlet boundary control given by
\begin{equation}
    \left\{
        \begin{aligned}
            y_t+Ay &= 0 \text{ in } Q_T, \\
            y &= Bu \text{ on } \Sigma_T, \\
            y(0)&=y_0 \text{ in } \Omega,
        \end{aligned}
    \right.
\end{equation}
where $B$ is defined as in the case of Neumann control just above.  In this case, \eqref{A1} can be verified by the method of transposition, and one arrives at the state space
\begin{equation}
    W=L^2(0,T;L^2(\Omega))\cap H^1(0,T;(H^1_0(\Omega)\cap H^2(\Omega))^*) \cap C([0,T];H^{-1}(\Omega)).
\end{equation}
This was carried out in, e.g., \cite[Thm.~2.1]{KV} with leading term in $A$ taken as the Laplacian for simplicity.

\begin{theorem}\label{thm:existence}
    There exists a minimizer $\bar u\in \calU$ to \eqref{eq:problem}.
\end{theorem}
\begin{proof}
    We first consider the case of $\calU=H^1(0,T;\R^2)$. Since $J_\beta$ is bounded from below, there exists a minimizing sequence $\{u_n\}_{n\in\N}$ that is bounded in $H^1(0,T;\R^2)$. Hence, by coercivity of $J_\beta$, there exists a subsequence, still denoted by $\{u_n\}_{n\in\N}$, with $u_n\rightharpoonup \bar u$ in $H^1(0,T;\R^2)$ and $u_n \to \bar u$ pointwise in $(0,T)$. This implies pointwise convergence of $|u_{n,1}(t)u_{n,2}(t)|\to |\bar u_1(t)\bar u_2(t)|$. Together with the continuity of $S$ and the weak lower semicontinuity of norms, this implies
    \begin{equation}
        J_\beta(\bar u ) \leq \liminf_{n\to\infty} J_\beta(u_n) = \inf_{u\in\calU} J_\beta(u),
    \end{equation}
    i.e., $\bar u$ is a minimizer.

    The case of $\calU$ finite dimensional follows similarly, since boundedness in $L^2(0,T;\R^N)$ then directly implies strong and hence pointwise convergence.
\end{proof}

We now address the convergence of solutions to \eqref{eq:problem} as $\beta\to\infty$ to a solution to the ``hard switching'' control problem
\begin{equation}\label{eq:problem_hard}
    \left\{   \begin{aligned}
            \min_{u \in \calU}\  &\frac12\norm{Su-y^d}_{L^2(0,T;L^2(\omobs))}^2 + \frac\alpha2 \norm{u}_{L^2(0,T;\R^2)}^2+ \frac\eps2 \norm{u_t}_{L^2(0,T;\R^2)}^2\\
            \text{s.t. }\ &u_1(t)u_2(t) = 0, \quad t\in [0,T].
    \end{aligned}\right.
\end{equation}
\begin{proposition}\label{thm:convergence_hard}
    The family $\{u_\beta\}_{\beta\geq 0}$ of minimizers to \eqref{eq:problem} contains at least one convergent sequence $\{u_{\beta_n}\}_{n\in\N}$ with $\beta_n\to\infty$.  
    The limit $\bar u\in\calU$ of every such sequence is a solution to \eqref{eq:problem_hard}.
\end{proposition}
\begin{proof}
    Again, we only consider the case of $\calU=H^1(0,T;\R^2)$, the other case being analogous. First, let $\{u_{\beta_n}\}_{n\in\N}$ with $\beta_n\to\infty$ be a sequence of minimizers to \eqref{eq:problem}. Since $J_\beta(u_\beta)< J_\beta(0)=J_0(0)$ for any $\beta>0$, this sequence is bounded in $H^1(0,T;\R^2)$ and hence contains a subsequence $\{u_n\}_{n\in\N}$, with $u_n\rightharpoonup \bar u$ in $H^1(0,T;\R^2)$ and $u_n \to \bar u$ pointwise in $(0,T)$.
    Furthermore, $J_\beta(u_\beta)< J_0(0)$ also implies 
    \begin{equation}
        \int_0^T |u_{n,1}(t)u_{n,2}(t)|\,dt \leq \beta_n^{-1} J_0(0) \to 0.
    \end{equation}
    Hence, 
    \begin{equation}
        \bar u_1(t)\bar u_2(t) = \lim_{n\to\infty} u_{n,1}(t)u_{n,2}(t) = 0 \qquad\text{for all }t\in[0,T].
    \end{equation}

    Now, let $\{u_n\}_{n\in\N}$ be any such sequence.
    Together with optimality of $u_n$, the above implies that for any $\tilde u\in\calU$ with $\tilde u_1(t)\tilde u_2(t)=0$ in $[0,T]$, we have
    \begin{equation}
        \begin{aligned}
            \frac12\norm{S\tilde u-y^d}_{L^2(0,T;L^2(\omobs))}^2 &+ \frac\alpha2 \norm{\tilde u}_{L^2(0,T;\R^2)}^2+ \frac\eps2 \norm{\tilde u_t}_{L^2(0,T;\R^2)}^2 \\[0.33em]
                                                                 &= J_{\beta_n}(\tilde u) \geq J_{\beta_n}(u_n)\geq J_0(u_n)
            \\&= \frac12\norm{S u_n-y^d}_{L^2(0,T;L^2(\omobs))}^2 + \frac\alpha2 \norm{u_n}_{L^2(0,T;\R^2)}^2+ \frac\eps2 \norm{u_{n,t}}_{L^2(0,T;\R^2)}^2.
        \end{aligned}
    \end{equation}
    Taking the limes inferior as $n\to\infty$ and using continuity of $S$ and weak lower semicontinuity of the norms now yields $J(\tilde u) \geq J(\bar u)$, i.e., $\bar u$ is a global minimizer.
\end{proof}

\section{Optimality conditions}\label{sec:optimality}

To derive optimality conditions, we can make use of the calculus of Clarke's generalized derivative \cite{Clarke:1990a,Clarke:2013}.
\begin{theorem}\label{thm:optimality}
    Any local minimizer $\bar u\in\calU$ to \eqref{eq:problem} satisfies
    \begin{equation}\label{eq:opt_incl}
    0\in S_0^*(S\bar u-z) + \alpha \bar u -\eps \bar u_{tt} + \beta\sign(\bar u_1\bar u_2)\begin{pmatrix} \bar u_2\\ \bar u_1 \end{pmatrix}.
    \end{equation}
\end{theorem}
\begin{proof}
    First, we consider the functional
    \begin{equation}
        G:L^2(0,T;\R^2)\to\R, \qquad u\mapsto \int_0^T|u_1(t)u_2(t)|\,dt.
    \end{equation}
    Since $H:L^1(0,T)\to\R$, $v\mapsto \int_0^T|v|\,dt$, is finite-valued, locally Lipschitz and convex, $H$ is regular at any $v\in L^1(0,T)$, and the generalized derivative coincides with the sub\-differential in the sense of convex analysis; see \cite[Prop.~2.2.7]{Clarke:1990a}.
    Furthermore, $T:L^2(0,T;\R^2)\to L^1(0,T)$, $u\mapsto u_1u_2$, is strictly differentiable.
    Hence, by \cite[Theorem 2.3.10]{Clarke:1990a}, $G=H\circ T$ is regular at any $u\in L^2(0,T;\R^2)$, and
    \begin{equation}
    \partial_C G(u) = T'(u)^*\partial H(T(u)) = \begin{pmatrix}u_2\\u_1\end{pmatrix}\sign(u_1u_2).
    \end{equation}
    Since $G$ is regular and the remaining terms in \eqref{eq:problem} are continuously Fréchet-differentiable, we can apply the sum rule for generalized gradients, e.g., from \cite[Prop.~2.3.3 with Cor.~1]{Clarke:1990a}, from which the desired result follows.
\end{proof}
Note that for $\calU=H^1(0,T;\R^2)$, the right-hand side of \eqref{eq:opt_incl} is to be understood as a subset of $H^1(0,T;\R^2)^*$.

\bigskip

Using directional derivatives, we can show that non-switching arcs can have a length of at most $\sqrt{\eps}$. In the following, we set $\norm{S_0}:=\norm{S_0}_{\mathcal{L}(L^2(0,T;\R^2),L^2(0,T;L^2(\omobs))}$ for brevity.
\begin{theorem}\label{thm:penalty_exact}
    If $\beta > (\norm{S_0}^2+\alpha+\pi^2)$, then $\bar u_1(t)\bar u_2(t) = 0$ for all $t\in[0,T]$ apart from intervals of length at most $\sqrt\eps$.
\end{theorem}
\begin{proof}
    Let $u\in L^2(0,T;\R^2)$ be given and
    assume that there exists $t_0\in (0,T)$ and $\delta >\sqrt\eps$ such that $u_1(t)u_2(t) \neq 0$ for all $t\in (t_0,t_0+\delta)$ and $u_1(t_0)u_2(t_0)=u_1(t_0+\delta)u_2(t_0+\delta)=0$. Without loss of generality, we can assume that both $u_1(t)>0$ and $u_2(t)>0$ for all $t\in (t_0,\delta)$.
    Furthermore, since $\beta > (\norm{S_0}^2+\alpha+\pi^2)$, there exists a $\rho\in(0,\delta/2)$ such that 
    \begin{equation}
        \beta > \left(\norm{S_0}^2+\alpha+\frac{\delta^2}{(\delta-2\rho)^2}\pi^2\right)>(\norm{S_0}^2+\alpha+\pi^2).
    \end{equation}
    Set
    \begin{equation}
        h_1(t) = \begin{cases}
            \sqrt{\frac2{\delta-2\rho}}\sin\left(\tfrac{\pi}{\delta-2\rho} (t-t_0-\rho)\right) & t\in I:=[t_0+\rho,t_0+\delta-\rho],\\
            0 &\text{else.}
        \end{cases}
    \end{equation}
    Then, $h_1 \in H^1(0,T)$ with
    \begin{equation}
        \norm{h_1}^2_{L^2(0,T)} = 1 \qquad\text{and}\qquad \norm{(h_1)_t}_{L^2(0,T)}^2 = \frac{\pi^2}{(\delta-2\rho)^2}.
    \end{equation}

    We now consider directional derivatives in the specific direction $h=(h_1,-h_1)$. For this purpose, we first introduce
    \begin{equation}
        g:\R\to\R,\qquad s \mapsto G(u+sh),
    \end{equation}
    and show that $g$ is differentiable in $0$. First, we have by definition of $h$ that
    \begin{equation}
        \begin{aligned}
            {g(s)-g(0)} &= \int_0^T \left[|(u_1+sh_1)(u_2-sh_1)| - |u_1u_2|\right]\,dt\\
                        &= \int_I \left[|(u_1+sh_1)(u_2-sh_1)| - |u_1u_2|\right]\,dt.
        \end{aligned}
    \end{equation}
    By the continuity of $u$ and $h$, there exists $\tilde s>0$ such that 
    \begin{equation}\label{eq:s_small}
        u_1(t) \pm sh_1(t) \geq 0 \quad \text{and}\quad u_2(t) \pm sh_1(t) \geq 0
        \qquad\text{for all }t\in I \text{ and } s\in[-\tilde s,\tilde s].
    \end{equation}
    Furthermore, $u_1(t)u_2(t)>0$ for all $t\in I$ by assumption.
    Hence,
    \begin{equation}\label{eq:ddiv1}
        g(s) - g(0) = \int_I \left[(u_1+sh_1)(u_2-sh_1) - u_1u_2\right]\,dt = \int_I \left[sh_1(u_2-u_1) -s^2 h_1^2\right]\,dt
    \end{equation}
    and therefore
    \begin{equation}\label{eq:ddiv2}
        g'(0) = \lim_{s\to 0} \frac{g(s) - g(0)}{s} = \int_I h_1(u_2-u_1)\,dt.
    \end{equation}
    For the second derivative $g''(0)$, we can proceed in the same way using \eqref{eq:s_small} to obtain
    \begin{equation}\label{eq:ddiv3}
        \begin{aligned}[t]
            g''(0) &= \lim_{s\to 0} \frac{g(s) -2 g(0) + g(-s)}{s^2} \\
                   &= \lim_{s\to 0} \frac1{s^2} \int_I \left[ (u_1+sh_1)(u_2-sh_1) - 2u_1u_2 + (u_1-sh_1)(u_2+sh_1)\right]\,dt\\
                   &= \lim_{s\to 0} \int_I -2h_1^2\,dt = -\norm{h}_{L^2(0,T;\R^2)}^2. 
        \end{aligned}
    \end{equation}
    Comparing \eqref{eq:ddiv1}, \eqref{eq:ddiv2}, and \eqref{eq:ddiv3}, we obtain that
    \begin{equation}\label{eq:dtaylor}
        g(s) = g(0) + sg'(0) + \frac{s^2}2 g''(0).
    \end{equation}    
    Since the remaining terms in the cost functional are differentiable, we have for 
    \begin{equation}
        j:\R\to\R,\qquad s\mapsto J(u+sh),
    \end{equation}
    that
    \begin{equation}
        \begin{aligned}
            j''(0) &= \norm{S_0h}_{L^2(0,T;L^2(\omobs))}^2 + \alpha\norm{h}^2_{L^2(0,T;\R^2)} + \eps \norm{h_t}^2_{L^2(0,T;\R^2)} - \beta \norm{h}^2_{L^2(0,T;\R^2)}\\
                   &\leq (\norm{S_0}^2 + \alpha - \beta)\norm{h}^2_{L^2(0,T;\R^2)} + \eps \norm{h_t}^2_{L^2(0,T;\R^2)}\\
                   &=2(\norm{S_0}^2 + \alpha - \beta) + 2\frac{\pi^2}{(\delta-2\rho)^2} \eps\\
                   &\leq 2\left(\norm{S_0}^2 + \alpha - \beta + \frac{\delta^2}{(\delta-2\rho)^2}\pi^2\right) <0
        \end{aligned}
    \end{equation}
    by the choice of $\delta$, $\rho$, and $\beta$.

    Now, from \eqref{eq:dtaylor} and the fact that the remaining terms in $J$ are quadratic and hence that the second-order Taylor expansion of $j$ is exact, it follows that for all $s$ with $sj'(0)\leq 0$, we have
    \begin{equation}
        \begin{aligned}
            J(u+sh) = j(s) = j(0) + sj'(0) + \frac{s^2}2 j''(0) < j(0) = J(u).
        \end{aligned}
    \end{equation}
    Hence, $u$ cannot be a local minimizer.
\end{proof}
The above proof relies on balancing the $L^2$ and $H^1$ norm for the perturbation $h$ to get an upper bound on the length of possible non-switching intervals.
Since the function $h_1$ used in this proof is the shifted and scaled first Dirichlet eigenfunction of the Laplacian on the interval $(0,1)$, which is the only function for which equality holds in the Poincaré inequality, the above result is likely to be sharp in this respect. 
Thus, perfect switching for $\eps>0$ cannot be guaranteed in general. 
However, it follows from Proposition~\ref{thm:convergence_hard} that non-switching arcs have to vanish for $\beta\to\infty$.
We also point out that the  condition on $\beta$ in \cref{thm:penalty_exact} is independent of the initial condition and possible further source terms.

We now turn to the finite-dimensional case, where $\calU$ consists of piecewise constant functions on a given grid
\begin{equation}
    0=t_1<\dots< t_M =T,
\end{equation}
and we can (and must) take $\eps=0$.
\begin{theorem}\label{thm:switching_pc}
    If $\,\calU$ consists of piecewise constant functions and $\beta > \norm{S_0}^2 + \alpha$, then $\bar u_1(t)\bar u_2(t) =0$ for all $t\in [0,T]$.
\end{theorem}
\begin{proof}
    We proceed as above. Let $u\in\calU$ with $u_1(t)u_2(t)\neq 0$ be given. Then there exists an interval $I_j:=(t_j,t_{j+1}]$ such that $u_1(t)u_2(t)\neq 0$ for all $t\in I_j$. As before, we can assume $u_1 >0$ and $u_2>0$ on $I_j$. We now choose 
    \begin{equation}
        h_1(t) = \begin{cases}
            1 & t\in I_j,\\
            0 &\text{else,}
        \end{cases}
    \end{equation}
    with $\norm{h_1}_{L^2(0,T)}^2 = \tau_j := t_{j+1}-t_j$ and consider again for $h=(h_1,-h_1)^T$ the function $g(s):= G(u+sh)$. We then obtain (using the fact that $u$ is constant on $I_j$) that
    \begin{equation}
        g(s) - g(0) = \tau_j \left[ |(u_1+s)(u_2-s)| - |u_1u_2|\right].
    \end{equation}
    Since $u_1$ and $u_2$ here are strictly positive constants, there again exists an $\tilde s$ such that $u_1\pm s>0$ and $u_2\pm s>0$ on for all $s\in [-\tilde s,\tilde s]$. Hence,
    \begin{equation}
        g'(0) = \tau_j (u_2-u_1)\qquad\text{and}\qquad g''(0) = -2\tau_j = -\norm{h}_{L ^2(0,T;\R^2)}^2.
    \end{equation}
    As above, the latter implies that
    \begin{equation}
        j''(0) \leq 2(\norm{S_0}^2 + \alpha - \beta) \tau_j < 0
    \end{equation}
    and hence that $J(u+sh) < J(u)$ for all $s$ with $sj'(0)\leq 0$.
\end{proof}

\begin{remark}\label{rem:switching_findim}
    The proof of \cref{thm:switching_pc} relies on the fact that for piecewise constant $u$, i.e., $u_i=\sum_{j=1}^M \xi^j_i \chi_{I_j}$ for $i\in\{1,2\}$, there holds
    \begin{equation}\label{eq:discrete_penalty}
        G(u) = \sum_{j=1}^M \tau_j |\xi^j_1\xi^j_2|.
    \end{equation}
    If $\calU$ is an arbitrary finite-dimensional subspace of $L^2(0,T;\R^2)$, i.e., $u_i=\sum_{j=1}^M \xi^j_i e_j$ for $i\in\{1,2\}$ and some basis functions $e_j$, and the switching penalty $G$ is replaced by the right-hand side of \eqref{eq:discrete_penalty}, the above proof can be modified to show that $\beta > C^{-1}(\norm{S_0}^2 + \alpha)$ implies that $\xi^j_1\xi^j_2=0$ for all $j$, where $C$ is the constant of equivalence for the discrete and continuous norm on $\calU\subset L^2(0,T;\R^2)$.

    The relation between pointwise switching and switching of the coefficients depends on the specific choice of $e_j$.
\end{remark}

\section{Numerical solution}\label{sec:solution}

The numerical solution is based on a primal-dual reformulation of the optimality condition \eqref{eq:opt_incl}, which states that there exists a
\begin{equation}\label{eq:def_q}
    \bar q \in \beta\sign(\bar u_1\bar u_2) = \partial (\beta\norm{\cdot}_{L^1})(\bar u_1\bar u_2)\subset L^\infty(0,1)
\end{equation}
such that
\begin{equation}\label{eq:opt_f}
S^*_0(S\bar u-z) + \alpha \bar u -\eps \bar u_{tt} + \bar q \begin{pmatrix}\bar u_2\\\bar u_1\end{pmatrix} = 0,
\end{equation}
where $\partial(\norm{\cdot}_{L^1})$ denotes the convex subdifferential of the $L^1(0,T)$ norm. Since we wish to formulate the algorithm in function space, we assume in this section that $\eps>0$ and hence $\calU=H^1(0,T;\R^2)$.

We now proceed as in the case of sparse control, e.g., \cite{CK:2009a,COAP:2011a,COCV:2015}:
By Fenchel duality applied to the convex functional $\norm{\cdot}_{L^1}$, \eqref{eq:def_q} is equivalent to
\begin{equation}\label{eq:dual_q}
    \bar u_1\bar u_2 \in \partial(\calI_{B_\infty(0,\beta)})(\bar q),
\end{equation}
where $\calI_{B_\infty(0,\beta)}$ denotes the indicator function (in the sense of convex analysis) of the closed ball in $L^\infty(0,T)$ around $0$ with radius $\beta$. The subdifferential on the right-hand side of \eqref{eq:dual_q} is now replaced by its Moreau--Yosida regularization in $L^2(0,T)$, which for any $\gamma>0$ is given by
\begin{equation}
    \partial(\calI_{B_\infty(0,\beta)})_\gamma(q) = \frac1\gamma \left(q- \prox_{\gamma\calI_{B_\infty(0,\beta)}}(q)\right) = \frac1\gamma\left(q-\proj_{B_\infty(0,\beta)}(q)\right).
\end{equation}
Here, $\prox_{f}$ denotes the proximal mapping of a convex function $f$, which for the indicator function of a convex set $C$ coincides with the metric projection $\proj_C$ onto $C$.
The Moreau--Yosida regularization of \eqref{eq:opt_incl} can thus be written as
\begin{equation}\label{eq:opt_gamma}
    \left\{\begin{aligned}
        S_0^*(Su_\gamma-z) + \alpha u_\gamma -\eps (u_\gamma)_{tt} + q_\gamma \begin{pmatrix}u_{\gamma,2}\\u_{\gamma,1}\end{pmatrix} &= 0\\[0.5em]
            \gamma u_{\gamma,1} u_{\gamma,2} - \max(0,q_\gamma-\beta) - \min(0,q_\gamma+\beta)&=0,
    \end{aligned}\right.
\end{equation}
where the second equation is to be understood pointwise almost everywhere in $(0,T)$.

It is known that the Moreau--Yosida regularization of the subdifferential of the indicator function is equivalent to adding an $L^2$ norm to the primal functional; see, e.g., \cite[Remark 3.2]{CK:2009a}. In our case, this leads to the problem
\begin{equation}\label{eq:problem_gamma}
    \begin{multlined}[t]
        \min_{u \in \calU} \frac12\norm{Su-y^d}_{L^2(0,T;L^2(\omobs))}^2 + \frac\alpha2 \norm{u}_{L^2(0,T;\R^2)}^2+ \frac\eps2 \norm{u_t}_{L^2(0,T;\R^2)}^2\\ + \beta \int_0^T |u_1(t)u_2(t)|\,dt + \frac\gamma2 \int_0^T |u_1(t)u_2(t)|^2\,dt .
    \end{multlined}
\end{equation}
A similar proof as for \cref{thm:existence} yields existence of a solution $u_\gamma$. Proceeding as in the derivation of \cref{thm:optimality}, we deduce the existence of
\begin{equation}\label{eq:gamma_inc}
    q_\gamma \in \beta\sign(u_{\gamma,1}u_{\gamma,2})+\gamma u_{\gamma,1}u_{\gamma,2} = \partial \left(\beta\norm{\cdot}_{L^1}+\frac\gamma2\norm{\cdot}_{L^2}^2\right)(u_{\gamma,1}u_{\gamma,2})
\end{equation}
such that the first relation of \eqref{eq:opt_gamma} holds. As shown in \cite[Remark 3.2]{CK:2009a}, the subdifferential inclusion \eqref{eq:gamma_inc} is equivalent to
\begin{equation}
    u_{\gamma,1}u_{\gamma,2} = \frac1\gamma \max(0,q_\gamma-\beta) + \frac1\gamma\min(0,q_\gamma+\beta), 
\end{equation}
which is a reformulation of the second relation of \eqref{eq:opt_gamma}. Hence, we obtain existence of a solution $(u_\gamma,q_\gamma)$
to \eqref{eq:opt_gamma}. A standard argument shows weak subsequential convergence of $u_\gamma$ as $\gamma \to 0$ to a minimizer $\bar u\in\calU$ to \eqref{eq:problem}.

\bigskip

Since $u_\gamma\in H^1(0,T;\R^2)$, it follows from \eqref{eq:gamma_inc} that $q_\gamma\in L^\infty(0,T)$. It is well-known (e.g., from \cite{Kunisch:2008a,Ulbrich:2011}) that the pointwise $\max$ and $\min$ are Newton differentiable from $L^p(0,T)$ to $L^1(0,T)$ for any $p>1$, with Newton derivative in direction $h\in L^p(0,T)$ given pointwise almost everywhere by
\begin{equation}
D_N \max(0,q-\beta)h = \begin{cases} h & \text{if }q\leq \beta,\\0 & \text{else,}\end{cases}
    \qquad
D_N \min(0,q+\beta)h = \begin{cases} h & \text{if }q\geq -\beta,\\0 & \text{else.}\end{cases}
\end{equation}
Thus, \eqref{eq:opt_gamma} considered as an operator equation from $H^1(0,T;\R^2)\times L^2(0,T)$ to $H^1(0,T;\R^2)^*\times L^2(0,T)$ is Newton differentiable. We therefore apply a semismooth Newton method for its solution.

\bigskip

In the numerical realization, we follow a homotopy approach. Since the problem is genuinely nonconvex, which is detrimental to the convergence behavior of Newton methods, this is more involved than in the convex case.
In particular, note that the Moreau--Yosida regularization acts as a smoothing of the indicator function; at the same time, it exacerbates the nonconvexity of the problem; cf. \eqref{eq:problem_gamma}. We therefore proceed as follows:
Starting from an initial guess $(u^0,q^0)=(0,0)$ and $\gamma=0$, we solve a sequence of problems for increasing $\gamma\geq 0$ and fixed $\beta=\beta_{\min}$ until a given final value $\gamma_{\max}>0$ is reached, taking the previous solution as an initial guess. Keeping this value of $\gamma=\gamma_{\max}$ fixed, we then similarly increase $\beta$ from some $\beta_{\min}$ to a value $\beta_{\max}$ which is chosen adaptively as the first $\beta$ such that the maximal pointwise violation $\swerror:=\|u_1u_2\|_\infty = \max_{t\in I} |u_1(t) u_2(t)|$ of the switching property drops below a given tolerance.
This adaptive choice is advantageous since it ensures large enough $\beta$ to obtain switching while avoiding too large $\beta$ that would complicate the numerical solution at no additional benefit.
Finally, keeping $\beta=\beta_{\max}>0$ fixed, $\gamma$ is reduced again until some specified $\gamma_{\min}$ is reached. 

In our experience, the homotopy parameters $\beta_{\min}$ and $\gamma_{\min}$ can be chosen very small in general, while $\gamma_{\max}$ has to be chosen above a critical value that depends on the problem. 
For any such choice, the obtained controls are then robust with respect to the value of $\gamma_{\max}$.

\section{Numerical examples}\label{sec:examples}
\renewcommand{\textfraction}{0.0}

We illustrate the properties of solutions to \eqref{eq:problem} for distributed control of the heat equation (i.e., $A=-\Delta$) on the unit square $\Omega=(-1,1)^2$ with homogeneous Neumann and initial conditions and different choices of control and observation configurations. The final time is always set to $T=10$. The spatial discretization uses standard piecewise linear finite elements, while the temporal discretization is also chosen as continuous and piecewise linear to obtain a conforming finite element approximation of the $H^1$ seminorm in \eqref{eq:problem} as well as the weak (temporal) Laplacian in \eqref{eq:opt_incl}.
Following \cref{rem:switching_findim}, the switching penalty is replaced by its discrete version (i.e., the right-hand side of \eqref{eq:discrete_penalty}), which allows a componentwise formulation of the second relation in \eqref{eq:gamma_inc}. We similarly replace the $L^2$ norm by its discrete version, which can be interpreted as a mass lumping.

The first example, based on the example from \cite{switchingcontrol}, illustrates the effect of the control cost parameter $\alpha$. In the second example, we investigate the influence of the $H^1$-regularization on the solution.
The MATLAB implementation of the proposed approach used to generate these results can be downloaded from \url{https://github.com/clason/nonconvexswitching}.

\subsection{Example 1}

To compare the proposed method with earlier work, we take up the example from \cite{switchingcontrol} with $N=2$ control components and $101$ equidistant time points.
For any $\alpha\geq 4\cdot 10^{-4}$, the convex solution method from \cite{switchingcontrol} yields switching controls,
which are unique global mininimizers to the respective problems. 
However, smaller values of $\alpha$ lead to non-switching solutions. In particular, we find for $\alpha=10^{-4}$ resp. $10^{-5}$ a total of $7$ resp. $40$ time points where both control components are essentially active.

In the following, we apply the proposed nonconvex approach for different values of the control cost parameter $\alpha$ with a fixed $\eps=10^{-7}$.
The homotopy parameters are set to $\beta_{\min}=10^{-5}$, $\gamma_{\min}=10^{-9}$, and $\gamma_{\max}=10^2$, and incremented or decremented by a factor of $10$. In this and the following examples, the final value $\beta_{\max}$ is chosen according to the relative tolerance $\swerror\leq 10^{-10} \max\{\|u_1\|_\infty,\|u_2\|_\infty\}$. 
For each pair $(\beta,\gamma)$, the maximum number of semismooth Newton iterations is set to $5$. While this leads to early termination of the semismooth Newton method in the first homotopy steps (to save numerical effort), we always observe convergence with relative tolerance $10^{-6}$ or absolute tolerance $10^{-7}$ in the residual norm of the optimality system \eqref{eq:opt_gamma} during the later steps of the homotopy. 

The solutions for different values of $\alpha$ are shown in \cref{fig:ex0_oc_diff_al}, where the solution for $\alpha=4\cdot 10^{-4}$ is very similar to the solution to the convex problem (differing only in the points $t_j$ with $u_1(t_j)=u_2(t_j)=0$, where in the convex case one of the control components is active).
Since the convex solution is computed with $\eps=0$, this indicates that
the influence of $\eps$ on the solution is negligible for this choice of $\alpha$ and $\eps$. For smaller values of $\alpha$, the solutions to the nonconvex problem maintain the switching property while increasing in amplitude until $\alpha \leq 10^{-8}$, after which any choice of $\alpha$ yields the same numerical solution.
\definecolor{mycolor1}{rgb}{0.231674, 0.318106, 0.544834}%
\definecolor{mycolor2}{rgb}{0.369214, 0.788888, 0.382914}%
\begin{figure}[p]
    \centering
    \begin{subfigure}[t]{0.495\textwidth}
        \begin{tikzpicture}

\begin{axis}[%
width=\textwidth,
xmin=0,
xmax=10,
ymin=-14,
ymax=7,
legend style={legend pos=south east,legend style={draw=none}},
]
\addplot [color=mycolor1,solid,line width=1.5pt]
  table[row sep=crcr]{%
0	-12.8947981725769\\
0.1	-12.9075534921147\\
0.2	-12.94452995823\\
0.3	-12.9960888367891\\
0.4	-13.048880241559\\
0.5	-13.0873519957121\\
0.6	-13.0929306584708\\
0.7	-13.0464640032147\\
0.8	-12.9284269553977\\
0.9	-12.7204627207735\\
1	-12.4006533207203\\
1.1	-11.7059693532323\\
1.2	-6.92475953403372e-06\\
1.3	-10.2013999243443\\
1.4	-3.44082660758588e-18\\
1.5	-5.36270133217392e-18\\
1.6	-5.87403864267641e-18\\
1.7	-4.90868022272044e-18\\
1.8	-5.13207600673656e-18\\
1.9	-4.58286342734824e-18\\
2	-3.96136258107628e-18\\
2.1	-4.26286060352836e-18\\
2.2	-8.66285117229395e-19\\
2.3	3.39090892141146e-17\\
2.4	2.18921239955677\\
2.5	3.15841654843332\\
2.6	3.94653554689907\\
2.7	4.56283280940429\\
2.8	4.99783691243395\\
2.9	5.24761977475226\\
3	5.31442197782043\\
3.1	5.20648677417569\\
3.2	4.9375482987951\\
3.3	4.5261203128017\\
3.4	3.99453403987568\\
3.5	3.36784712604461\\
3.6	2.67261364767134\\
3.7	1.9356481102388\\
3.8	1.18289466529244\\
3.9	0.444113153369306\\
4	3.48235409273157e-16\\
4.1	6.57847063659789e-17\\
4.2	9.21290260453136e-17\\
4.3	-2.0161503078578\\
4.4	-2.595192234292e-06\\
4.5	-2.73263075484985e-06\\
4.6	-3.29295325509682\\
4.7	-3.67014488680082\\
4.8	-3.94879408928233\\
4.9	-4.21042486148728\\
5	-4.4727669534275\\
5.1	-4.7515464354185\\
5.2	-5.0608607553431\\
5.3	-5.41219629711001\\
5.4	-5.81351757071123\\
5.5	-6.26858523175173\\
5.6	-6.77648209161216\\
5.7	-7.33141653626513\\
5.8	-7.92278293733574\\
5.9	-8.53550473669519\\
6	-9.15062153761945\\
6.1	-9.74611409657754\\
6.2	-10.2979049197588\\
6.3	-10.7810078912065\\
6.4	-11.1707446485051\\
6.5	-11.443994185855\\
6.6	-11.5803813545553\\
6.7	-11.563379894046\\
6.8	-11.3812314503004\\
6.9	-11.0276839903998\\
7	-10.5024486337274\\
7.1	-9.81143049932533\\
7.2	-8.96661332931322\\
7.3	-7.98575263886642\\
7.4	-6.89162919609116\\
7.5	-5.70969253035583\\
7.6	-4.39986028797239\\
7.7	-3.50847606213271e-17\\
7.8	-5.07832695741617e-17\\
7.9	-1.24892116621736e-16\\
8	0.419844999538029\\
8.1	1.45033114234091\\
8.2	2.37615459755556\\
8.3	3.16315091605506\\
8.4	3.79851643618627\\
8.5	4.27560015262022\\
8.6	4.59356529972098\\
8.7	4.757153180585\\
8.8	4.7760959120637\\
8.9	4.66447617502588\\
9	4.4396033248979\\
9.1	4.12118753661403\\
9.2	3.72975926562953\\
9.3	3.2860867717283\\
9.4	2.80900019829269\\
9.5	2.31600953743379\\
9.6	1.81863946770737\\
9.7	1.31065732808323\\
9.8	-1.16499768936437e-16\\
9.9	-4.15765611538527e-17\\
10	-8.64770101805556e-17\\
};
\addlegendentry{$u_1$};

\addplot [color=mycolor2,solid,line width=1.5pt]
  table[row sep=crcr]{%
0	2.2930541876001e-18\\
0.1	6.66823127836299e-19\\
0.2	1.04689851844177e-18\\
0.3	-8.47227515067352e-20\\
0.4	-5.47765381071869e-20\\
0.5	4.01493787846483e-19\\
0.6	-1.26144057795981e-18\\
0.7	-1.21755880876589e-18\\
0.8	-3.84058952863221e-18\\
0.9	-4.22711804954289e-18\\
1	2.31012592093142e-18\\
1.1	-8.39632364331797e-19\\
1.2	-7.21603808772096e-06\\
1.3	2.34257872445874e-18\\
1.4	-9.99634361206388\\
1.5	-9.48621556913856\\
1.6	-8.64191127681026\\
1.7	-7.71150121264296\\
1.8	-6.7139092115807\\
1.9	-5.66650083403556\\
2	-4.5895076765538\\
2.1	-3.50553084706967\\
2.2	-2.43822210810337\\
2.3	-1.40180766602777\\
2.4	3.21979347331485e-18\\
2.5	2.0309490407676e-18\\
2.6	2.29884568981653e-18\\
2.7	2.98668534914459e-18\\
2.8	2.41649671481439e-18\\
2.9	2.15993180660747e-18\\
3	2.32439712514435e-18\\
3.1	2.85595937680733e-18\\
3.2	4.17637848608349e-18\\
3.3	5.59018382506872e-18\\
3.4	5.5560359663592e-18\\
3.5	9.63202726943149e-18\\
3.6	9.69691702573426e-18\\
3.7	2.01589466478454e-17\\
3.8	6.68948902966359e-17\\
3.9	1.09856227148936e-16\\
4	-0.595549383857826\\
4.1	-1.14923376460346\\
4.2	-1.62273402407012\\
4.3	1.50432407650058e-16\\
4.4	-2.65216134796303e-06\\
4.5	-2.84386523167689e-06\\
4.6	3.43354918615529e-17\\
4.7	3.82607602470184e-17\\
4.8	3.63492522608075e-18\\
4.9	8.9142463280147e-18\\
5	1.04551649663755e-17\\
5.1	1.01225342632128e-17\\
5.2	9.88629976700927e-18\\
5.3	1.04444617045592e-17\\
5.4	9.06468138388265e-18\\
5.5	9.40330728194609e-18\\
5.6	4.86835421063454e-18\\
5.7	6.86400555388647e-18\\
5.8	2.9224915937416e-18\\
5.9	4.02065578881505e-18\\
6	2.25306855989934e-18\\
6.1	1.49659251761048e-18\\
6.2	6.40383257979049e-19\\
6.3	4.75033346647399e-19\\
6.4	3.81279461436756e-20\\
6.5	5.27825432274307e-19\\
6.6	1.38614061206506e-19\\
6.7	4.53887773865959e-19\\
6.8	-3.28072329944169e-19\\
6.9	-5.72086006026371e-19\\
7	-1.55674584981438e-18\\
7.1	-2.37246106541703e-18\\
7.2	-3.51881809739403e-18\\
7.3	-4.90850984230255e-18\\
7.4	-7.10913419518615e-18\\
7.5	-1.02721845976838e-17\\
7.6	-2.20636423961484e-17\\
7.7	-3.3514571509772\\
7.8	-2.42145748027189\\
7.9	-1.39468885468885\\
8	2.15587891545483e-15\\
8.1	-8.09035389841139e-17\\
8.2	-4.27523817479759e-17\\
8.3	-2.75018367125452e-17\\
8.4	-2.03402688991139e-17\\
8.5	-1.71990628435072e-17\\
8.6	-1.57339069123165e-17\\
8.7	-1.44080379048017e-17\\
8.8	-1.350838415272e-17\\
8.9	-1.35570580288463e-17\\
9	-1.27160428693796e-17\\
9.1	-1.5929914116831e-17\\
9.2	-1.59675665177585e-17\\
9.3	-2.03833822524424e-17\\
9.4	-2.22364296614824e-17\\
9.5	-2.70056089304158e-17\\
9.6	-3.42750020034139e-17\\
9.7	-5.85102721560099e-17\\
9.8	0.826938472157433\\
9.9	0.397871784727592\\
10	0.182755077335884\\
};
\addlegendentry{$u_2$}
\end{axis}
\end{tikzpicture}%
        \caption{$\alpha=4\cdot10^{-4}$}
    \end{subfigure}
    \hfill
    \begin{subfigure}[t]{0.495\textwidth}
        \raggedleft
        \begin{tikzpicture}

\begin{axis}[%
width=\textwidth,
xmin=0,
xmax=10,
ymin=-30,
ymax=30,
legend style={legend pos=south east,legend style={draw=none}},
]
\addplot [color=mycolor1,solid,line width=1.5pt]
  table[row sep=crcr]{%
0	-24.5651321961312\\
0.1	-24.619085444937\\
0.2	-24.8116894207771\\
0.3	-25.1373490343833\\
0.4	-25.544862494599\\
0.5	-25.9755223405552\\
0.6	-26.3577229285311\\
0.7	-26.6200012644621\\
0.8	-26.6880284435546\\
0.9	-26.4839825736889\\
1	-25.8220944071374\\
1.1	-23.1343889557167\\
1.2	2.40008307205794e-16\\
1.3	1.92927578923713e-16\\
1.4	2.03086331460434e-16\\
1.5	2.21198444302313e-16\\
1.6	2.494330841044e-16\\
1.7	2.94974375100392e-16\\
1.8	3.65959211754588e-16\\
1.9	4.77135229446511e-16\\
2	2.23191616886065e-16\\
2.1	6.4253522194715\\
2.2	10.5268158840004\\
2.3	14.1829278409107\\
2.4	17.4627726147389\\
2.5	20.2943687996108\\
2.6	22.6078034712004\\
2.7	24.3499295100973\\
2.8	25.4878424810938\\
2.9	26.0100976284949\\
3	25.9271561852238\\
3.1	25.2705156468533\\
3.2	24.091052354796\\
3.3	22.4562987985227\\
3.4	20.4470906299833\\
3.5	18.1535209396532\\
3.6	15.670602100902\\
3.7	13.093693662569\\
3.8	10.5140634867056\\
3.9	8.01468797840399\\
4	5.66659919071184\\
4.1	3.52585883671563\\
4.2	1.63119032040697\\
4.3	-1.9637417023691e-10\\
4.4	-1.35439087775215\\
4.5	-2.45896392689039\\
4.6	-3.34234194953825\\
4.7	-4.0497138110112\\
4.8	-4.63683381888127\\
4.9	-5.16601292525641\\
5	-5.70196607019445\\
5.1	-6.30757493370328\\
5.2	-7.03974709828316\\
5.3	-7.94561912320926\\
5.4	-9.05927551999243\\
5.5	-10.3991897794086\\
5.6	-11.9664939459454\\
5.7	-13.7442064745174\\
5.8	-15.6974286917619\\
5.9	-17.7745466723567\\
6	-19.9093395069762\\
6.1	-22.0239423259342\\
6.2	-24.0324615388831\\
6.3	-25.8451283014292\\
6.4	-27.3727057974458\\
6.5	-28.5310192046444\\
6.6	-29.2452714625009\\
6.7	-29.4540571875389\\
6.8	-29.1127087144826\\
6.9	-28.1960063791138\\
7	-26.6998502438878\\
7.1	-24.6421555099978\\
7.2	-22.0624257770515\\
7.3	-19.0206282568632\\
7.4	-15.593059720625\\
7.5	-11.8519381327682\\
7.6	-7.64640323917124\\
7.7	7.17319964871731e-15\\
7.8	7.14000177091716e-15\\
7.9	3.77490697940223\\
8	7.32672439547615\\
8.1	10.5341865949825\\
8.2	13.3176687000108\\
8.3	15.6167633184388\\
8.4	17.3909512470077\\
8.5	18.6202834509273\\
8.6	19.3054847475085\\
8.7	19.4669845989028\\
8.8	19.1428210440266\\
8.9	18.3861481747824\\
9	17.2612542878117\\
9.1	15.8403143152616\\
9.2	14.1979056119429\\
9.3	12.4082286713072\\
9.4	10.5361505982429\\
9.5	8.61541626537377\\
9.6	6.36900759551184\\
9.7	5.57065832024707e-06\\
9.8	2.64673090574948e-14\\
9.9	7.94498617545268e-14\\
10	-1.50082759243566e-13\\
};
\addlegendentry{$u_1$};
\addplot [color=mycolor2,solid,line width=1.5pt]
  table[row sep=crcr]{%
0	3.1216185489666e-16\\
0.1	2.4695196618736e-16\\
0.2	2.38670866652221e-16\\
0.3	2.32624267667952e-16\\
0.4	2.21279691324311e-16\\
0.5	2.06267802896632e-16\\
0.6	1.93289786009582e-16\\
0.7	1.85073148199591e-16\\
0.8	1.76814777625785e-16\\
0.9	1.72014133891407e-16\\
1	1.70804301409846e-16\\
1.1	2.21013860418827e-16\\
1.2	-21.7366788135657\\
1.3	-22.1741508221229\\
1.4	-20.5773770536196\\
1.5	-18.5311806771306\\
1.6	-16.1736396950737\\
1.7	-13.5389942006553\\
1.8	-10.665868715311\\
1.9	-7.60029633869803\\
2	-4.32759460639503\\
2.1	2.32151963636537e-16\\
2.2	9.20880256492155e-17\\
2.3	2.08281013018023e-17\\
2.4	-1.08066582779088e-17\\
2.5	-2.69208237746416e-17\\
2.6	-3.80114114427523e-17\\
2.7	-4.75307550504141e-17\\
2.8	-5.73012579324631e-17\\
2.9	-6.79272055555667e-17\\
3	-8.10355199765365e-17\\
3.1	-9.89868246716596e-17\\
3.2	-1.23554055649966e-16\\
3.3	-1.58241939377967e-16\\
3.4	-2.09511749386957e-16\\
3.5	-2.87585027777156e-16\\
3.6	-4.14894521908069e-16\\
3.7	-6.34163574537965e-16\\
3.8	-1.03954950682675e-15\\
3.9	-1.87765967962457e-15\\
4	-3.91603080792967e-15\\
4.1	-1.04602900352194e-14\\
4.2	-5.12840629553421e-14\\
4.3	0.0198143057450128\\
4.4	-7.63048941021493e-14\\
4.5	-2.31372092531211e-14\\
4.6	-1.25684212710735e-14\\
4.7	-8.49918164241907e-15\\
4.8	-6.39802470531871e-15\\
4.9	-5.03303766931882e-15\\
5	-3.98067749333967e-15\\
5.1	-3.10566783461571e-15\\
5.2	-2.35599718796736e-15\\
5.3	-1.73646563663518e-15\\
5.4	-1.2399397648493e-15\\
5.5	-8.65247013742363e-16\\
5.6	-5.96480835148581e-16\\
5.7	-4.09504054111683e-16\\
5.8	-2.80483510357151e-16\\
5.9	-1.90994582493983e-16\\
6	-1.28630603324662e-16\\
6.1	-8.50178439601277e-17\\
6.2	-5.46050294020366e-17\\
6.3	-3.17686119829255e-17\\
6.4	-1.43856974428001e-17\\
6.5	3.99009000028032e-19\\
6.6	1.34980639021583e-17\\
6.7	2.66769573464819e-17\\
6.8	4.15316080747298e-17\\
6.9	5.94584416272238e-17\\
7	8.24680732281853e-17\\
7.1	1.15336582268976e-16\\
7.2	1.66872335163733e-16\\
7.3	2.56621926524582e-16\\
7.4	4.32272324458325e-16\\
7.5	8.58500909890557e-16\\
7.6	2.50588784478642e-15\\
7.7	-4.4232346487918\\
7.8	-1.89408761208171\\
7.9	5.30042150609118e-15\\
8	2.05285245464373e-15\\
8.1	1.05227678684353e-15\\
8.2	6.81889323973694e-16\\
8.3	5.09489613881633e-16\\
8.4	4.19474529426845e-16\\
8.5	3.68323383569379e-16\\
8.6	3.29511755838782e-16\\
8.7	2.85439353351682e-16\\
8.8	2.48785616860393e-16\\
8.9	2.18294152964647e-16\\
9	1.8687811953902e-16\\
9.1	1.29595664127263e-16\\
9.2	9.13472665954222e-18\\
9.3	-2.15044775541962e-16\\
9.4	-5.59277489298375e-16\\
9.5	2.61812849411795e-15\\
9.6	5.37059670722539e-15\\
9.7	5.60511253008263e-06\\
9.8	2.9482587740477\\
9.9	1.59765221881805\\
10	0.827328530353417\\
};
\addlegendentry{$u_2$};
\end{axis}
\end{tikzpicture}%
        \caption{$\alpha=1\cdot10^{-4}$}
    \end{subfigure}

    \bigskip 

    \begin{subfigure}[t]{0.495\textwidth}
        \begin{tikzpicture}

\begin{axis}[%
width=\textwidth,
xmin=0,
xmax=10,
ymin=-130,
ymax=130,
legend style={legend pos=south east,legend style={draw=none}},
]
\addplot [color=mycolor1,solid,line width=1.5pt]
  table[row sep=crcr]{%
0	-14.2464517111559\\
0.1	-15.2266446334971\\
0.2	-18.3906408548097\\
0.3	-23.9092042827919\\
0.4	-31.6056210425911\\
0.5	-41.1017715645197\\
0.6	-51.8316035356686\\
0.7	-63.004963933302\\
0.8	-73.4231292161766\\
0.9	-81.2011908474551\\
1	-83.3386594145327\\
1.1	-75.2559662752313\\
1.2	-50.3809012632528\\
1.3	3.08939551079368e-20\\
1.4	1.22323527869413e-24\\
1.5	-3.37498106924157e-24\\
1.6	2.08293671151958e-25\\
1.7	-2.327670843131e-24\\
1.8	7.42277908725037e-25\\
1.9	9.12264369253682e-25\\
2	-1.33568989412e-20\\
2.1	31.4713605182639\\
2.2	56.4127801799689\\
2.3	76.0386946674698\\
2.4	91.36981380249\\
2.5	103.034310397864\\
2.6	111.325757832421\\
2.7	116.325441455563\\
2.8	118.030094665132\\
2.9	116.44942408576\\
3	111.675393911299\\
3.1	103.917930139243\\
3.2	93.5199731632984\\
3.3	80.9547269594377\\
3.4	66.8161660789946\\
3.5	51.8049084221526\\
3.6	36.7162746544703\\
3.7	22.4246912250042\\
3.8	9.86479940048262\\
3.9	4.23696137256197e-18\\
4	-12.9335186054467\\
4.1	-22.3835305428879\\
4.2	-28.5646428136822\\
4.3	-31.7402668904453\\
4.4	-32.2260117479951\\
4.5	-30.40055635336\\
4.6	-26.71352548554\\
4.7	-21.6870904961253\\
4.8	-15.9028777231014\\
4.9	-9.97392698301471\\
5	-4.49860471574918\\
5.1	1.33032122724674e-16\\
5.2	2.05450241437278e-24\\
5.3	-1.57660256819561e-24\\
5.4	-1.97728930460764e-24\\
5.5	-2.20737839816281e-24\\
5.6	-1.69857757207221e-24\\
5.7	-2.1917851916139e-25\\
5.8	-8.49473130061332e-21\\
5.9	-20.8870496170449\\
6	-37.7551912567466\\
6.1	-52.2568436359346\\
6.2	-65.3703223362954\\
6.3	-77.4472145991699\\
6.4	-88.3848793287863\\
6.5	-97.8071844487124\\
6.6	-105.218804160603\\
6.7	-110.119252952993\\
6.8	-112.083476713851\\
6.9	-110.811643607885\\
7	-106.158068151714\\
7.1	-98.142909026085\\
7.2	-86.9528232946132\\
7.3	-72.9323697412046\\
7.4	-56.5693271397665\\
7.5	-38.4757437279633\\
7.6	-19.3626395114376\\
7.7	-1.49449479901072e-24\\
7.8	19.4483339777213\\
7.9	37.7168886323196\\
8	54.2182712902058\\
8.1	68.4468716307448\\
8.2	80.0010454883979\\
8.3	88.6094829917679\\
8.4	94.1478133982614\\
8.5	96.6373701221171\\
8.6	96.2180768587736\\
8.7	93.0884363806951\\
8.8	87.3991348820716\\
8.9	79.0981329881535\\
9	67.7213231554018\\
9.1	52.1710357508576\\
9.2	30.5311518830523\\
9.3	-1.01748137068667e-18\\
9.4	-1.4764282335666e-24\\
9.5	1.18881685122715e-24\\
9.6	-6.69908959989199e-24\\
9.7	2.70040957884611e-24\\
9.8	4.79659528849424e-25\\
9.9	4.30164927047011e-24\\
10	-4.81497100598269e-24\\
};
\addlegendentry{$u_1$};
\addplot [color=mycolor2,solid,line width=1.5pt]
  table[row sep=crcr]{%
0	-6.7505519647601e-25\\
0.1	2.87249947958652e-24\\
0.2	-2.19156387318175e-24\\
0.3	-1.8042048423332e-24\\
0.4	-2.0809409242924e-24\\
0.5	-1.70359653290241e-24\\
0.6	-2.46154407388922e-24\\
0.7	-2.45379044965569e-24\\
0.8	-2.06283427810355e-24\\
0.9	1.44288457389242e-24\\
1	-5.47334172210531e-25\\
1.1	-6.81981248419841e-26\\
1.2	3.91928513985281e-20\\
1.3	-48.8623390825396\\
1.4	-72.6180500157139\\
1.5	-80.700112881786\\
1.6	-79.6772980120803\\
1.7	-73.4984311305971\\
1.8	-63.8150704769679\\
1.9	-50.2177473427262\\
2	-30.3442074683575\\
2.1	1.90561733616904e-21\\
2.2	-3.17865478873323e-25\\
2.3	9.4911594179524e-26\\
2.4	2.85084029560657e-24\\
2.5	-1.60212727391943e-24\\
2.6	-1.89946678392496e-24\\
2.7	-2.4379918488117e-26\\
2.8	1.32251989045462e-25\\
2.9	3.79980249732605e-25\\
3	1.90588846344907e-25\\
3.1	-1.89242868615222e-25\\
3.2	-1.5199377063007e-24\\
3.3	8.35555802786232e-25\\
3.4	9.98927023338781e-25\\
3.5	1.55311436696376e-24\\
3.6	-6.82693064122734e-25\\
3.7	-3.68618476111394e-24\\
3.8	4.76432664978239e-19\\
3.9	-4.8509255641754\\
4	-2.13387782620305e-20\\
4.1	2.68429398563567e-24\\
4.2	-1.95535551662919e-24\\
4.3	1.38811242141379e-24\\
4.4	1.14217271134098e-24\\
4.5	7.51588973497407e-25\\
4.6	2.54216422050351e-24\\
4.7	-7.74286579819025e-25\\
4.8	-1.65660440655683e-24\\
4.9	-3.03827556817535e-24\\
5	3.01821633604019e-24\\
5.1	2.15035224188948\\
5.2	7.35902340508141\\
5.3	13.5976201766764\\
5.4	19.5046962148549\\
5.5	24.0539098531704\\
5.6	26.1988804516209\\
5.7	24.5060508382938\\
5.8	16.8066831706022\\
5.9	-1.36931478602028e-19\\
6	2.9879498232075e-24\\
6.1	-2.21063593005716e-25\\
6.2	-1.12285648975491e-24\\
6.3	3.68073513098927e-25\\
6.4	-9.75516186218529e-25\\
6.5	1.50253413339961e-24\\
6.6	1.26752914383593e-25\\
6.7	3.93609811046937e-24\\
6.8	1.85960578266964e-24\\
6.9	-5.64656552714865e-25\\
7	1.05534457598122e-24\\
7.1	-2.53787194336663e-25\\
7.2	-7.81913718084356e-25\\
7.3	-1.80287858352614e-24\\
7.4	1.39314016057456e-24\\
7.5	-3.87119231510426e-24\\
7.6	2.63932653950357e-21\\
7.7	-4.42065745865692\\
7.8	-1.00503284253138e-24\\
7.9	-3.0868544515067e-24\\
8	1.15539275923625e-24\\
8.1	-4.53602504791633e-24\\
8.2	-1.01481815351282e-25\\
8.3	-8.83201449101861e-25\\
8.4	3.48170367753589e-24\\
8.5	-1.3203743715762e-24\\
8.6	8.11970127657951e-25\\
8.7	9.93275299667731e-26\\
8.8	-6.85448196141837e-25\\
8.9	-1.52485836412369e-25\\
9	1.3438596693874e-24\\
9.1	1.64427370234358e-24\\
9.2	-7.5823097259367e-20\\
9.3	23.1721798612641\\
9.4	29.9185449730772\\
9.5	27.933335528333\\
9.6	22.4102772160505\\
9.7	16.5335133255954\\
9.8	12.0087259215029\\
9.9	9.50647666759477\\
10	8.79139720764568\\
};
\addlegendentry{$u_2$};
\end{axis}
\end{tikzpicture}%
        \caption{$\alpha=1\cdot10^{-6}$}
    \end{subfigure}
    \hfill
    \begin{subfigure}[t]{0.495\textwidth}
        \begin{tikzpicture}

\begin{axis}[%
width=\textwidth,
xmin=0,
xmax=10,
ymin=-130,
ymax=130,
legend style={legend pos=south east,legend style={draw=none}},
]
\addplot [color=mycolor1,solid,line width=1.5pt]
  table[row sep=crcr]{%
0	2.59374310193048e-25\\
0.1	1.82956635065091e-19\\
0.2	-10.8337490491751\\
0.3	-21.0582026124363\\
0.4	-31.4930420467097\\
0.5	-42.6149674451668\\
0.6	-54.4655896205691\\
0.7	-66.5487742147921\\
0.8	-77.6684509849852\\
0.9	-85.7000071172258\\
1	-87.3010197440223\\
1.1	-77.7168267187243\\
1.2	-50.869960256441\\
1.3	9.18170479933861e-21\\
1.4	-3.89159778415999e-25\\
1.5	1.37774292165626e-24\\
1.6	2.43148278747431e-24\\
1.7	-5.8299741801962e-25\\
1.8	-1.5812597748426e-24\\
1.9	4.09007230199563e-24\\
2	4.61773576957728e-21\\
2.1	32.6688222078998\\
2.2	59.1429552364506\\
2.3	80.0604696792038\\
2.4	96.2755614389049\\
2.5	108.441453731026\\
2.6	116.935134872052\\
2.7	121.913651633461\\
2.8	123.419372146577\\
2.9	121.479466233129\\
3	116.187038453195\\
3.1	107.75094652581\\
3.2	96.5238589190877\\
3.3	83.0109403737064\\
3.4	67.8705103432984\\
3.5	51.9077255037707\\
3.6	36.0654170347922\\
3.7	21.3978874890553\\
3.8	9.01865038678342\\
3.9	1.62877189654836e-17\\
4	-14.9125883685333\\
4.1	-25.6151147854697\\
4.2	-32.3730838722603\\
4.3	-35.5775654554089\\
4.4	-35.6845735799815\\
4.5	-33.1994041901454\\
4.6	-28.6828753255081\\
4.7	-22.7664355460137\\
4.8	-16.1578169778284\\
4.9	-9.63079765612532\\
5	-3.98898246204908\\
5.1	1.65422199115628e-24\\
5.2	5.41373212161295e-24\\
5.3	-4.06980796201298e-24\\
5.4	4.72779880084849e-25\\
5.5	-1.95518791636551e-24\\
5.6	1.27674064440131e-24\\
5.7	3.50642188512423e-24\\
5.8	3.77856103814223e-21\\
5.9	-19.9960833408252\\
6	-37.1404435571879\\
6.1	-52.3838012266863\\
6.2	-66.3718512387382\\
6.3	-79.317427231525\\
6.4	-91.0643894603241\\
6.5	-101.210119894573\\
6.6	-109.234500583076\\
6.7	-114.608185182906\\
6.8	-116.876094011123\\
6.9	-115.712748635695\\
7	-110.956355621956\\
7.1	-102.624350645975\\
7.2	-90.9166496403494\\
7.3	-76.2085542138436\\
7.4	-59.0364728933269\\
7.5	-40.0772569754381\\
7.6	-20.1172646880248\\
7.7	-6.57042980582332e-25\\
7.8	20.6930017923787\\
7.9	40.019747021723\\
8	57.3781932034291\\
8.1	72.2615089921612\\
8.2	84.2723924161386\\
8.3	93.146345696866\\
8.4	98.7662457573569\\
8.5	101.157451666522\\
8.6	100.454630151875\\
8.7	96.8353941961659\\
8.8	90.4141140860085\\
8.9	81.1062152054086\\
9	68.4777939431197\\
9.1	51.6474944219954\\
9.2	29.3143982380021\\
9.3	-1.04543469871698e-18\\
9.4	8.79413770034037e-25\\
9.5	5.10450158701134e-24\\
9.6	-2.02821728220079e-24\\
9.7	6.08519185684294e-24\\
9.8	4.33576680875429e-24\\
9.9	3.26741593561006e-24\\
10	6.07009415303138e-24\\
};
\addlegendentry{$u_1$};
\addplot [color=mycolor2,solid,line width=1.5pt]
  table[row sep=crcr]{%
0	-10.9909187448049\\
0.1	-8.42694968755426\\
0.2	1.21251112038228e-18\\
0.3	-8.04360792610484e-24\\
0.4	2.08473346091666e-25\\
0.5	-3.68312712100911e-24\\
0.6	-2.51418961707511e-24\\
0.7	7.23738041156684e-25\\
0.8	5.09879019753709e-25\\
0.9	-2.48736278042025e-24\\
1	-2.57496566415445e-24\\
1.1	2.38204937948985e-24\\
1.2	-8.42007825033041e-21\\
1.3	-48.7323912977117\\
1.4	-74.5586093114248\\
1.5	-84.7302832612769\\
1.6	-85.0954957992039\\
1.7	-79.5047239281779\\
1.8	-69.6269948614873\\
1.9	-54.9806776634687\\
2	-33.1026025589656\\
2.1	3.75664075190958e-20\\
2.2	1.43764404966589e-24\\
2.3	6.9149125542209e-25\\
2.4	2.03735349060328e-24\\
2.5	5.55131251866864e-25\\
2.6	2.36190568413732e-24\\
2.7	1.10923006270454e-24\\
2.8	4.47672832870911e-25\\
2.9	-1.58645051540662e-24\\
3	1.72202373023067e-24\\
3.1	1.26878306019959e-24\\
3.2	-1.78701251814471e-24\\
3.3	-1.96802779299576e-24\\
3.4	-1.96800067724267e-24\\
3.5	-9.71008687381712e-25\\
3.6	-4.8520156216227e-25\\
3.7	-1.60824186819633e-25\\
3.8	7.80582385361588e-18\\
3.9	-5.7704333450734\\
4	-9.19102712484208e-19\\
4.1	-6.40215313293552e-26\\
4.2	1.54486410514354e-25\\
4.3	2.3628929248976e-24\\
4.4	-2.24039975490395e-24\\
4.5	2.77107105761863e-24\\
4.6	-7.76277941885122e-25\\
4.7	-1.42050564152091e-24\\
4.8	2.18029643921294e-24\\
4.9	4.80545596675318e-25\\
5	1.77599443083502e-24\\
5.1	2.3683725601619\\
5.2	8.2778133780323\\
5.3	15.6003604409523\\
5.4	22.7592717498561\\
5.5	28.4429320928382\\
5.6	31.242762492664\\
5.7	29.2882481411992\\
5.8	19.9763661921536\\
5.9	-5.89668934548227e-19\\
6	2.18913077304674e-24\\
6.1	-1.5077658054836e-24\\
6.2	9.43492443608546e-25\\
6.3	1.64973556566501e-24\\
6.4	3.1244659482608e-24\\
6.5	-4.03033233693418e-24\\
6.6	-1.55828621588358e-24\\
6.7	-3.39816361805148e-24\\
6.8	7.66203445580174e-25\\
6.9	-1.35538833502949e-24\\
7	-1.63325965373836e-24\\
7.1	-3.42729378678626e-24\\
7.2	1.9101259954179e-24\\
7.3	-1.244391239394e-24\\
7.4	-4.0562769668842e-25\\
7.5	3.42814039517818e-24\\
7.6	-9.86448237369332e-25\\
7.7	-4.49359160383416\\
7.8	2.8041499478101e-24\\
7.9	4.39481353875389e-25\\
8	3.92545412198518e-24\\
8.1	-3.0305914933587e-25\\
8.2	-4.71165610170718e-25\\
8.3	-5.90861701409037e-26\\
8.4	1.49338411477093e-24\\
8.5	4.20654488665823e-24\\
8.6	3.03926691820423e-24\\
8.7	-2.48263645321083e-24\\
8.8	2.82150746356599e-24\\
8.9	-1.60599031422714e-24\\
9	6.22255838895104e-27\\
9.1	2.03322406805319e-25\\
9.2	-7.62240730613021e-20\\
9.3	22.2364666908659\\
9.4	29.3170302511121\\
9.5	27.7833504292506\\
9.6	22.5955414980488\\
9.7	17.0566785983077\\
9.8	13.0062212605176\\
9.9	11.0565723549253\\
10	10.6122708533356\\
};
\addlegendentry{$u_2$};
\end{axis}
\end{tikzpicture}%
        \caption{$\alpha=1\cdot10^{-8}$}
    \end{subfigure}
    \caption{Optimal switching controls for different values of $\alpha$}
    \label{fig:ex0_oc_diff_al}
\end{figure}
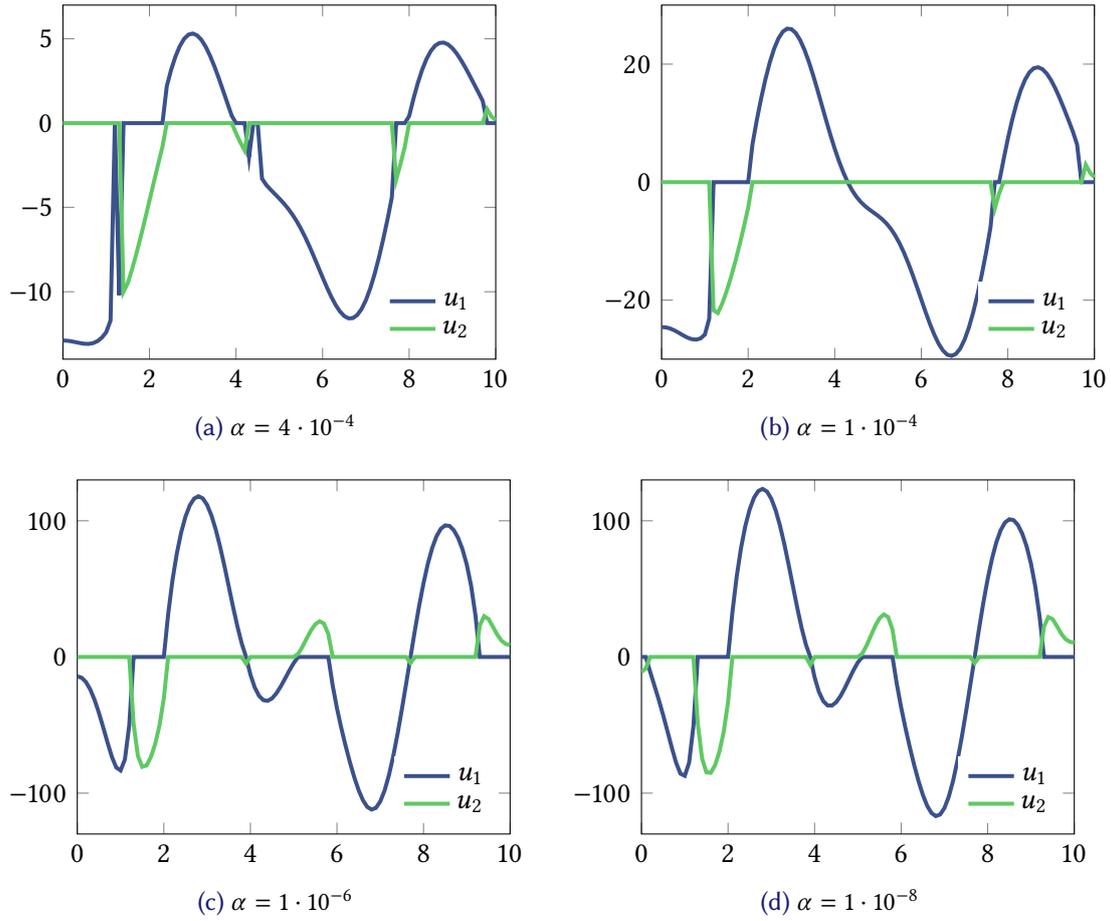
This is illustrated in more detail in \cref{tab:results_ex0}, which shows in particular that the tracking error (and therefore the optimal functional value $\bar{J}$) can be reduced significantly by decreasing $\alpha$. At the same time, the number $N_{\mathrm{sw}}$ of ``switching points'' (i.e., points $t_j$ such that $|u_1(t_{j})| \geq |u_2(t_{j})|$ but $|u_1(t_{j+1})|<|u_2(t_{j+1})|$ or vice versa)
increases only slightly due to the presence of the $H^1$-regularization.
(This relation will be investigated more closely in the next example.)
The fourth column gives the switching error $\swerror=\max_{t\in I}|u_1 (t)u_2(t)|$, confirming quantitatively that the switching property is satisfied up to a high accuracy for all cases.
The last two columns address the convergence of the combined homotopy and semismooth Newton method by giving the residual norm of the optimality system \eqref{eq:opt_gamma} for the final switching penalty $\beta=\beta_{\max}$ and $\gamma=10^{-9}$. Note that in all cases, we have $\beta_{\max}>\alpha$ and hence a genuinely nonconvex problem.
\begin{table}[p]
    \caption{Results for different values of $\alpha$ (Example 1): optimal value $\bar{J}$, number of switching points $N_{\mathrm{sw}}$, switching error $\swerror$, residual norm of the optimality system \eqref{eq:opt_gamma}, final switching penalty $\beta_{\max}$}\label{tab:results_ex0}
    \centering
    \begin{tabular}{llclll}
        \toprule
        $\alpha$ & $\bar{J}$ & $N_{\mathrm{sw}}$ & $\swerror$ & optimality & $\beta_{\max}$\\
        \midrule
        $4\cdot 10^{-4}$ & $0.762$ & $11$  & $5\cdot 10^{-11}$ & $7\cdot 10^{-15}$ & $10^3$\\  
        $10^{-4}$        & $0.581$ & $7$  & $3\cdot 10^{-11}$ & $8\cdot 10^{-14}$ & $10^2$\\
        $10^{-5}$        & $0.209$ & $9$  & $4\cdot 10^{-11}$ & $4\cdot 10^{-14}$ & $10^2$\\
        $10^{-6}$        & $0.069$ & $9$  & $3\cdot 10^{-16}$ & $1\cdot10^{-14}$        & $10^{-4}$\\
        $10^{-7}$        & $0.050$ & $10$ & $2\cdot 10^{-10}$ & $2\cdot 10^{-13}$ & $10^1$\\
        $10^{-8}$        & $0.048$ & $10$ & $7\cdot 10^{-17}$ & $2\cdot 10^{-14}$ & $10^{-4}$\\
        $10^{-9}$        & $0.048$ & $10$ & $7\cdot 10^{-17}$ & $2\cdot 10^{-14}$ & $10^{-4}$\\
        \bottomrule
    \end{tabular}
\end{table}

Finally, the convergence history for a run of the semismooth Newton method for $\alpha=10^{-8}$ and two representative pairs $(\beta,\gamma)$ (corresponding to low and high penalization, respectively) is shown in \cref{tab:results_ex0_superlin}. The superlinear decay of the norm $\|F(u^k,q^k)\|$ of the residual in \eqref{eq:opt_gamma} in each iteration $k$ can be observed clearly.
\begin{table}[t]
    \caption{Convergence history of the semismooth Newton method (Example 1) for $\alpha=10^{-8}$: residual norm $\|F(u^k,q^k)\|$ over iteration $k$}\label{tab:results_ex0_superlin}
    \begin{subtable}[t]{\textwidth}
        \centering
        \caption{$\beta=10^{-5}$, $\gamma=10^{-5}$}
        \begin{tabular}{ccccccc}
            \toprule
            $k$         & $0$                  & $1$                  & $2$                  & $3$                  & $4$ & $5$\\
            \midrule
            $\|F(u^k,q^k)\|$          & $2.7\cdot 10^{-5}$ & $9.8\cdot 10^{-6}$ & $3.6\cdot 10^{-6}$ & $1.1\cdot 10^{-7}$ & $8.0\cdot 10^{-10}$ &\phantom{$4.4\cdot10^{-10}$}\\
            \bottomrule
        \end{tabular}
    \end{subtable}

    \bigskip 

    \begin{subtable}[t]{\textwidth}
        \centering
        \caption{$\beta=10^{-1}$, $\gamma=10^{2}$}
        \begin{tabular}{ccccccc}
            \toprule
            $k$         & $0$                  & $1$                  & $2$                  & $3$                  & $4$ & $5$\\
            \midrule
            $\|F(u^k,q^k)\|$          & $1.5\cdot 10^{-3}$ & $3.7\cdot 10^{-4}$ & $9.2\cdot 10^{-5}$ & $4.8\cdot 10^{-5}$ & $4.8\cdot 10^{-6}$ & $4.4\cdot 10^{-10}$\\
            \bottomrule
        \end{tabular}
    \end{subtable}
\end{table}

\subsection{Example 2}

\begin{figure}[p]
    \centering
    \begin{tikzpicture}

\begin{axis}[%
width=0.4\textwidth,
scale only axis,
xmin=0,
xmax=10,
ymin=0,
ymax=22,
legend style={legend pos=north west,legend style={draw=none}},
]
\addplot [color=mycolor1,solid,line width=1.5pt]
  table[row sep=crcr]{%
0	0\\
0.1	0.000310889676799488\\
0.2	0.00493508882334778\\
0.3	0.0246567096710148\\
0.4	0.0765001870088557\\
0.5	0.182372542187894\\
0.6	0.367290024609168\\
0.7	0.657306816048791\\
0.8	1.07728382129735\\
0.9	1.64864710997755\\
1	2.38728757031316\\
1.1	3.30174563942213\\
1.2	4.39180805142067\\
1.3	5.64761844200694\\
1.4	7.04937193113662\\
1.5	8.56762745781211\\
1.6	10.164232941279\\
1.7	11.7938197208773\\
1.8	13.4057866110226\\
1.9	14.9466625730374\\
2	16.3627124296868\\
2.1	17.6026337878861\\
2.2	18.620186427436\\
2.3	19.376598311396\\
2.4	19.8426049159664\\
2.5	20\\
2.6	19.8426049159664\\
2.7	19.376598311396\\
2.8	18.620186427436\\
2.9	17.6026337878861\\
3	16.3627124296868\\
3.1	14.9466625730374\\
3.2	13.4057866110226\\
3.3	11.7938197208773\\
3.4	10.164232941279\\
3.5	8.56762745781211\\
3.6	7.04937193113662\\
3.7	5.64761844200694\\
3.8	4.39180805142068\\
3.9	3.30174563942213\\
4	2.38728757031316\\
4.1	1.64864710997755\\
4.2	1.07728382129735\\
4.3	0.657306816048793\\
4.4	0.367290024609167\\
4.5	0.182372542187895\\
4.6	0.0765001870088562\\
4.7	0.0246567096710147\\
4.8	0.00493508882334782\\
4.9	0.000310889676799483\\
5	4.49855881137988e-63\\
5.1	0.000310889676799479\\
5.2	0.00493508882334778\\
5.3	0.0246567096710146\\
5.4	0.0765001870088554\\
5.5	0.182372542187893\\
5.6	0.367290024609168\\
5.7	0.657306816048791\\
5.8	1.07728382129735\\
5.9	1.64864710997755\\
6	2.38728757031316\\
6.1	3.30174563942212\\
6.2	4.39180805142067\\
6.3	5.64761844200693\\
6.4	7.04937193113662\\
6.5	8.5676274578121\\
6.6	10.164232941279\\
6.7	11.7938197208773\\
6.8	13.4057866110226\\
6.9	14.9466625730374\\
7	16.3627124296868\\
7.1	17.6026337878861\\
7.2	18.620186427436\\
7.3	19.376598311396\\
7.4	19.8426049159664\\
7.5	20\\
7.6	19.8426049159664\\
7.7	19.376598311396\\
7.8	18.620186427436\\
7.9	17.6026337878861\\
8	16.3627124296868\\
8.1	14.9466625730374\\
8.2	13.4057866110226\\
8.3	11.7938197208773\\
8.4	10.164232941279\\
8.5	8.56762745781211\\
8.6	7.04937193113663\\
8.7	5.64761844200696\\
8.8	4.39180805142066\\
8.9	3.30174563942212\\
9	2.38728757031316\\
9.1	1.64864710997756\\
9.2	1.07728382129736\\
9.3	0.657306816048789\\
9.4	0.367290024609168\\
9.5	0.182372542187895\\
9.6	0.0765001870088564\\
9.7	0.0246567096710153\\
9.8	0.0049350888233477\\
9.9	0.000310889676799486\\
};
\addlegendentry{$u_1$};
\addplot [color=mycolor2,solid,line width=1.5pt]
  table[row sep=crcr]{%
0	10\\
0.1	9.96137346387737\\
0.2	9.84623899471337\\
0.3	9.65681246565741\\
0.4	9.39672242703368\\
0.5	9.07091524951596\\
0.6	8.68552626289954\\
0.7	8.24772068423965\\
0.8	7.76550901155084\\
0.9	7.24754229914006\\
1	6.70289330551129\\
1.1	6.14082989724222\\
1.2	5.57058729051458\\
1.3	5.00114571027522\\
1.4	4.44101984565615\\
1.5	3.89806608576945\\
1.6	3.37931294473099\\
1.7	2.89081934671727\\
1.8	2.43756456397911\\
1.9	2.02337261026362\\
2	1.65087281971106\\
2.1	1.32149722017001\\
2.2	1.03551417465352\\
2.3	0.792096650382418\\
2.4	0.589422415885803\\
2.5	0.42480249556895\\
2.6	0.29483335787942\\
2.7	0.195567603897884\\
2.8	0.122697379567451\\
2.9	0.0717443734231883\\
3	0.0382500935044278\\
3.1	0.0179601470466927\\
3.2	0.00699647333607276\\
3.3	0.002011896443561\\
3.4	0.00032195724373366\\
3.5	9.73450173947918e-06\\
3.6	2.49341021956449e-07\\
3.7	0.000102038383179388\\
3.8	0.00101454728610358\\
3.9	0.00430110496579768\\
4	0.0123283548355074\\
4.1	0.0281741080655795\\
4.2	0.0555066111390811\\
4.3	0.0984391537330873\\
4.4	0.161364753970941\\
4.5	0.248776320944915\\
4.6	0.365078188268672\\
4.7	0.514395221686577\\
4.8	0.7003858185631\\
4.9	0.926065033547981\\
5	1.19364378515658\\
5.1	1.50438963079513\\
5.2	1.8585139570942\\
5.3	2.25508963787452\\
5.4	2.69200228802391\\
5.5	3.16593721641073\\
5.6	3.67240308625403\\
5.7	4.20579216084964\\
5.8	4.7594758810974\\
5.9	5.32593342382961\\
6	5.89690986043864\\
6.1	6.46359960561065\\
6.2	7.01685004491669\\
6.3	7.5473795824793\\
6.4	8.04600387610155\\
6.5	8.50386374196641\\
6.6	8.91264812334556\\
6.7	9.2648056306839\\
6.8	9.55373847076875\\
6.9	9.77397308119721\\
7	9.92130245798318\\
7.1	9.99289598852871\\
7.2	9.98737355425671\\
7.3	9.90484171796586\\
7.4	9.74689092839085\\
7.5	9.51655382444445\\
7.6	9.21822586909536\\
7.7	8.857550652785\\
7.8	8.44127324485279\\
7.9	7.97706590696275\\
8	7.4733312865187\\
8.1	6.93898885607941\\
8.2	6.38325083720013\\
8.3	5.81539412969394\\
8.4	5.24453485162826\\
8.5	4.67941197914663\\
8.6	4.1281862622972\\
8.7	3.5982600933647\\
8.8	3.09612333338419\\
8.9	2.62722928145916\\
9	2.19590402571034\\
9.1	1.80529137346291\\
9.2	1.45733445386445\\
9.3	1.15279395265416\\
9.4	0.891301811280405\\
9.5	0.671448135767677\\
9.6	0.490898048156811\\
9.7	0.346534306143932\\
9.8	0.234620742598557\\
9.9	0.150980959675279\\
};
\addlegendentry{$u_2$};
\end{axis}
\end{tikzpicture}%
    \caption{Generating control $u_d$ for Example 2}\label{fig:generating_control}
\end{figure}
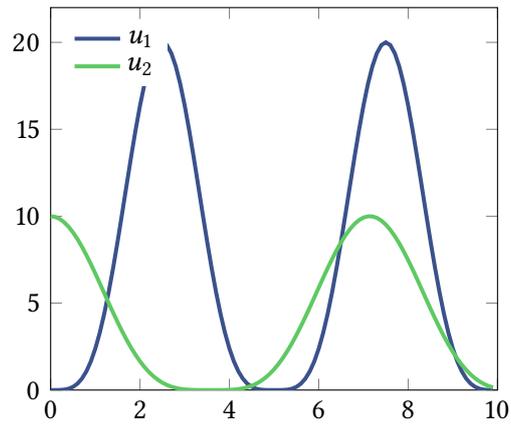

The second example addresses the influence of the $H^1$-regularization parameter $\eps$ on the solution.
Here, we set
\begin{equation}
    \omega_1 = \{ (x_1,x_2)\in \Omega: x_1\leq 0 \}, \qquad \omega_2 = \{(x_1,x_2)\in\Omega: x_1>0\},
\end{equation}
and $Bu=(\chi_{\omega_1}(x)u_1(t)+\chi_{\omega_2}(x)u_2(t))/10$. 
The desired state $y^d$ is the solution corresponding to the control 
\begin{equation}
    u_d = \left(20\sin^4(2\pi t/T),\,10\cos^4(1.4\pi t/T)\right),
\end{equation}
see \cref{fig:generating_control}. By construction, this desired state is expected to be difficult to attain by a pure switching control since both controls
contribute significantly during the second half of the time interval.
To investigate the influence of $\eps$ on the solution, we fix $\alpha=10^{-6}$, which corresponds to comparatively small control costs.
The homotopy parameters are set to $\beta_{\min}=10^{-5}$, $\gamma_{\min}=10^{-9}$, and $\gamma_{\max}=10^4$, and are incremented or decremented by a factor of $10$.
The maximal number of semismooth Newton iterations is again set to $5$; we observe the same convergence behavior as in the first example.

The solutions for different values of $\eps$ are shown in \cref{fig:oc_diff_eta}.
We see that all solutions exhibit switching, which is consistent with \cref{thm:penalty_exact}.
During the first half of the time interval, we observe for all $\eps>0$ the switching behavior expected from the generating control $u_d$.
During the second half, the number of switching points increases for $\eps\to 0$.
\begin{figure}[p]
    \centering
    \begin{subfigure}[t]{0.495\textwidth}
        \begin{tikzpicture}

\begin{axis}[%
width=\textwidth,
xmin=0,
xmax=10,
ymin=-5,
ymax=25,
legend style={legend pos=north west,legend style={draw=none}},
]
\addplot [color=mycolor1,solid,line width=1.5pt]
  table[row sep=crcr]{%
0	3.30872245021211e-24\\
0.1	2.41260508948977e-23\\
0.2	1.46938195186239e-23\\
0.3	1.41185304360881e-23\\
0.4	2.18222443567069e-23\\
0.5	1.68489161255227e-23\\
0.6	1.25604579096641e-23\\
0.7	9.251051940898e-24\\
0.8	2.2606684905555e-23\\
0.9	-9.67949229815457e-15\\
1	2.79687112050896\\
1.1	5.32635575621676\\
1.2	7.60203749083186\\
1.3	9.63910539471043\\
1.4	11.4525068334414\\
1.5	13.0556494563384\\
1.6	14.4596125454499\\
1.7	15.6728015513691\\
1.8	16.700976891124\\
1.9	17.5475833747659\\
2	18.2143051164577\\
2.1	18.7017708917689\\
2.2	19.0103374732706\\
2.3	19.1408834161605\\
2.4	19.0955528217199\\
2.5	18.8783978357972\\
2.6	18.4958793512197\\
2.7	17.9571974460528\\
2.8	17.2744356284039\\
2.9	16.4625156883409\\
3	15.5389719444595\\
3.1	14.5235646263263\\
3.2	13.4377612253605\\
3.3	12.3041217034311\\
3.4	11.1456278940977\\
3.5	9.98499932001566\\
3.6	8.84403673668282\\
3.7	7.74303133064996\\
3.8	6.70027174055883\\
3.9	5.73167352169506\\
4	4.8505466932674\\
4.1	4.06750734373145\\
4.2	3.39052941280724\\
4.3	2.82512348943866\\
4.4	2.37462123014069\\
4.5	2.04053744953108\\
4.6	1.82297736144172\\
4.7	1.72105423633173\\
4.8	1.7332829141455\\
4.9	1.8579172325797\\
5	2.09320422190192\\
5.1	2.43753463214583\\
5.2	2.88947746454391\\
5.3	3.44769520716761\\
5.4	4.11074575085887\\
5.5	4.87678593257646\\
5.6	5.74319965392568\\
5.7	6.70618007349792\\
5.8	7.76029996657752\\
5.9	8.89810671658044\\
6	10.1097783148548\\
6.1	11.3828742257429\\
6.2	12.702210105759\\
6.3	14.0498784795602\\
6.4	15.4054289411921\\
6.5	16.7462118327073\\
6.6	18.0478792079906\\
6.7	19.2850268909002\\
6.8	20.4319521743567\\
6.9	21.4634938053691\\
7	22.3559148303133\\
7.1	23.08778506488\\
7.2	23.6408186261099\\
7.3	24.0006232556837\\
7.4	24.1573219819327\\
7.5	24.106013827193\\
7.6	23.847048371899\\
7.7	23.386098581292\\
7.8	22.7340267727398\\
7.9	21.9065493433447\\
8	20.9237161999317\\
8.1	19.8092301326664\\
8.2	18.5896390431773\\
8.3	17.2934395350357\\
8.4	15.9501335141636\\
8.5	14.5892799666136\\
8.6	13.2395819021401\\
8.7	11.928043747006\\
8.8	10.6792274711254\\
8.9	9.51462693167901\\
9	8.45216977373027\\
9.1	7.50584547200592\\
9.2	6.68544726799148\\
9.3	5.99640571539515\\
9.4	5.4396826664326\\
9.5	5.01168780516216\\
9.6	4.70417503934236\\
9.7	4.50407469468001\\
9.8	4.39321757273873\\
9.9	4.3479137523608\\
10	4.33832821987245\\
};
\addlegendentry{$u_1$};
\addplot [color=mycolor2,solid,line width=1.5pt]
  table[row sep=crcr]{%
0	13.2745896472445\\
0.1	13.1578609372075\\
0.2	12.8059939364948\\
0.3	12.2131032641074\\
0.4	11.3684731209876\\
0.5	10.2571118365584\\
0.6	8.86054115648517\\
0.7	7.15785884022421\\
0.8	5.12713727303017\\
0.9	2.74722396918092\\
1	-5.6676051820815e-15\\
1.1	-5.59430255762988e-23\\
1.2	-3.77422007273453e-24\\
1.3	6.72856438834005e-24\\
1.4	1.87698138330114e-23\\
1.5	2.84771856687933e-26\\
1.6	-1.0333642007925e-23\\
1.7	-2.60425064497316e-24\\
1.8	1.32171563046713e-23\\
1.9	1.63762803425963e-23\\
2	6.0632757338594e-24\\
2.1	-1.04029978639497e-23\\
2.2	-1.73286458199002e-23\\
2.3	-1.45798614879644e-23\\
2.4	-1.39936191061771e-23\\
2.5	1.32074310742121e-24\\
2.6	2.57339164223175e-24\\
2.7	7.11593918800892e-24\\
2.8	-1.1075431900969e-24\\
2.9	9.71882889522752e-25\\
3	-1.2208038012555e-23\\
3.1	3.46284294227365e-24\\
3.2	8.79449210081527e-24\\
3.3	8.71211456664664e-24\\
3.4	1.91090160850512e-23\\
3.5	2.54781032455948e-23\\
3.6	5.27816332128515e-24\\
3.7	1.83575947050845e-23\\
3.8	3.21838911076544e-23\\
3.9	1.52125313183598e-23\\
4	6.17304542899809e-23\\
4.1	6.87254987009887e-23\\
4.2	9.77391708896877e-23\\
4.3	1.2966225981207e-22\\
4.4	1.40894494912231e-22\\
4.5	1.72704289702068e-22\\
4.6	2.13431771541182e-22\\
4.7	1.88773957210338e-22\\
4.8	1.59985976718815e-22\\
4.9	1.30227415317643e-22\\
5	6.62818365250947e-23\\
5.1	5.46554721802457e-23\\
5.2	2.62921142231269e-24\\
5.3	-1.26117693021042e-23\\
5.4	-2.56416798591946e-23\\
5.5	-1.16715315724114e-23\\
5.6	8.50604384206973e-24\\
5.7	-2.6806694489255e-23\\
5.8	-1.54603601954171e-23\\
5.9	-1.30511436636545e-23\\
6	-3.86303120094335e-24\\
6.1	-3.63499954997139e-23\\
6.2	-6.02752010168979e-24\\
6.3	-1.07754007867063e-23\\
6.4	-1.44057931207087e-23\\
6.5	-7.53137819125233e-24\\
6.6	-5.265423075327e-24\\
6.7	-1.12852879063263e-23\\
6.8	-1.49963527169316e-24\\
6.9	1.67921188276194e-24\\
7	-9.86390144718635e-24\\
7.1	1.25085162029608e-24\\
7.2	2.97448946728268e-24\\
7.3	-5.6145802556269e-24\\
7.4	-2.2294173544983e-23\\
7.5	-1.4522463331352e-23\\
7.6	2.34762163653972e-24\\
7.7	4.54308670267215e-25\\
7.8	-1.24225529793626e-23\\
7.9	-5.04165728346158e-24\\
8	-6.54683358342471e-24\\
8.1	-1.37452609364455e-23\\
8.2	-6.41432057605563e-24\\
8.3	2.09141485298947e-24\\
8.4	-1.13463108652986e-23\\
8.5	-1.05551025295261e-23\\
8.6	-9.8355367891164e-24\\
8.7	3.99992577075379e-25\\
8.8	-1.39647674528677e-23\\
8.9	-1.43112093777436e-23\\
9	-2.32537068019434e-23\\
9.1	-9.08203771243369e-24\\
9.2	-2.76197528072199e-23\\
9.3	-2.63926577036833e-23\\
9.4	-2.97483546180544e-23\\
9.5	-2.34598249494962e-23\\
9.6	-3.96826604817593e-23\\
9.7	-1.82886246375009e-23\\
9.8	-1.0978622367216e-23\\
9.9	-1.39639917509142e-24\\
10	-8.70502351396313e-24\\
};
\addlegendentry{$u_2$};
\end{axis}
\end{tikzpicture}%
        \caption{$\eps=10^{-3}$}
    \end{subfigure}
    \hfill
    \begin{subfigure}[t]{0.495\textwidth}
        \begin{tikzpicture}

\begin{axis}[%
width=\textwidth,
xmin=0,
xmax=10,
ymin=-5,
ymax=30,
legend style={legend pos=north west,legend style={draw=none}},
]
\addplot [color=mycolor1,solid,line width=1.5pt]
  table[row sep=crcr]{%
0	1.26123133390314e-16\\
0.1	6.26921755518217e-17\\
0.2	6.16174561998969e-17\\
0.3	5.84275859900268e-17\\
0.4	5.41453685360417e-17\\
0.5	5.36474707790933e-17\\
0.6	5.30841449701531e-17\\
0.7	5.27141155815561e-17\\
0.8	4.6221727760674e-17\\
0.9	5.00668502398682e-17\\
1	8.27077023662194e-17\\
1.1	8.0158257453756e-16\\
1.2	3.91718101031717\\
1.3	7.14155645863339\\
1.4	9.77844736010155\\
1.5	11.9422478536929\\
1.6	13.7372093938815\\
1.7	15.246987557133\\
1.8	16.530790370752\\
1.9	17.624002338015\\
2	18.5415946641921\\
2.1	19.2830010036306\\
2.2	19.8374977208931\\
2.3	20.1894160069656\\
2.4	20.322752686885\\
2.5	20.2249256772909\\
2.6	19.8895561418539\\
2.7	19.3182556089705\\
2.8	18.5214651462906\\
2.9	17.5184400656591\\
3	16.3365049242515\\
3.1	15.0097228252614\\
3.2	13.5771335844371\\
3.3	12.0807183231879\\
3.4	10.5632446319156\\
3.5	9.06613684496022\\
3.6	7.62750070855757\\
3.7	6.28041111147346\\
3.8	5.05154654845808\\
3.9	3.96022548081246\\
4	3.01786967127464\\
4.1	2.22788931076102\\
4.2	1.58595739333698\\
4.3	1.08061761807393\\
4.4	0.694155565483029\\
4.5	0.403656371547941\\
4.6	0.182181828509133\\
4.7	2.84852043288774e-14\\
4.8	5.82531349419816e-16\\
4.9	-1.09656464344049e-16\\
5	-3.9241280457358e-16\\
5.1	-4.18565468856136e-16\\
5.2	-3.2101621697363e-16\\
5.3	-1.79393828648715e-16\\
5.4	-8.59199170328264e-17\\
5.5	-4.55311302719298e-17\\
5.6	-2.93411286436715e-17\\
5.7	-1.57937525149441e-17\\
5.8	-4.63812909619294e-18\\
5.9	-8.36146530773872e-19\\
6	-5.77622855438454e-18\\
6.1	-9.99273642752693e-18\\
6.2	-5.06277821590924e-18\\
6.3	-1.61727431526246e-16\\
6.4	6.74859523362257\\
6.5	12.1710147166579\\
6.6	16.4457287762441\\
6.7	19.770454074513\\
6.8	22.3301913719404\\
6.9	24.2801930314609\\
7	25.739849218844\\
7.1	26.793722825203\\
7.2	27.4968063468269\\
7.3	27.8817791918428\\
7.4	27.9667101120631\\
7.5	27.7621712148519\\
7.6	27.2771518537439\\
7.7	26.5234626252222\\
7.8	25.5185378789631\\
7.9	24.2866852266335\\
8	22.8589206264712\\
8.1	21.2715779898261\\
8.2	19.5639129385729\\
8.3	17.7749405130917\\
8.4	15.9397703325754\\
8.5	14.0857374991602\\
8.6	12.2286831524099\\
8.7	10.369820448568\\
8.8	8.4937337169988\\
8.9	6.56820032848491\\
9	4.54668684496959\\
9.1	2.37453981344028\\
9.2	2.40212763203987e-05\\
9.3	2.01419001441604e-15\\
9.4	7.43314882585128e-16\\
9.5	3.65693154755678e-17\\
9.6	-7.62144609611037e-17\\
9.7	-9.11522413398656e-17\\
9.8	-5.3612848135162e-17\\
9.9	-3.21267397660309e-17\\
10	-3.50225635121741e-17\\
};
\addlegendentry{$u_1$};
\addplot [color=mycolor2,solid,line width=1.5pt]
  table[row sep=crcr]{%
0	11.2149871737416\\
0.1	11.2176839567804\\
0.2	11.2287024255437\\
0.3	11.2455991473562\\
0.4	11.2482399707857\\
0.5	11.1959552445254\\
0.6	11.0235466103113\\
0.7	10.6368685672777\\
0.8	9.90917778450564\\
0.9	8.67992572902587\\
1	6.75824777383052\\
1.1	3.93406121752233\\
1.2	-4.39824444374417e-16\\
1.3	3.39655746760665e-17\\
1.4	3.69596337196388e-18\\
1.5	-6.37316710062454e-18\\
1.6	-9.77755362843277e-18\\
1.7	-1.11470531004134e-17\\
1.8	-1.26142350371918e-17\\
1.9	-1.37007676816366e-17\\
2	-1.32486609204733e-17\\
2.1	-1.25094089189888e-17\\
2.2	-1.24007083454962e-17\\
2.3	-1.31688543843818e-17\\
2.4	-1.36353788031714e-17\\
2.5	-1.35497400670607e-17\\
2.6	-1.34711250925847e-17\\
2.7	-1.39671702934777e-17\\
2.8	-1.49627574206973e-17\\
2.9	-1.5565329223367e-17\\
3	-1.6481418968953e-17\\
3.1	-1.86072998233908e-17\\
3.2	-2.07219393376741e-17\\
3.3	-2.26238824169993e-17\\
3.4	-2.48594612679315e-17\\
3.5	-2.90293656190814e-17\\
3.6	-3.49712828453407e-17\\
3.7	-4.26426560755266e-17\\
3.8	-5.31612739099412e-17\\
3.9	-6.72698044753902e-17\\
4	-8.50039335223743e-17\\
4.1	-1.34509989692651e-16\\
4.2	-3.57581405042746e-16\\
4.3	-1.06827081406172e-15\\
4.4	-2.81507025477249e-15\\
4.5	-7.61304632078866e-15\\
4.6	-1.60943885657889e-13\\
4.7	0.292879346764633\\
4.8	0.643696010110395\\
4.9	1.0672279685441\\
5	1.5821795271204\\
5.1	2.20917553203897\\
5.2	2.96805352347822\\
5.3	3.87410724183041\\
5.4	4.93290861806374\\
5.5	6.13338233772074\\
5.6	7.43892291189698\\
5.7	8.77657221741731\\
5.8	10.0246355383068\\
5.9	10.9996670470662\\
6	11.4445072771253\\
6.1	11.0200826240719\\
6.2	9.30490078642172\\
6.3	5.80772653558842\\
6.4	2.02023222605931e-16\\
6.5	2.99143429459144e-17\\
6.6	2.89758370423003e-17\\
6.7	2.76267128705714e-17\\
6.8	2.65316459031003e-17\\
6.9	2.42290984379485e-17\\
7	2.00130654839255e-17\\
7.1	1.2964875783219e-17\\
7.2	1.09103804874308e-18\\
7.3	-1.89518974513313e-17\\
7.4	-1.72190594352589e-17\\
7.5	1.6644397833204e-17\\
7.6	-8.68293044155803e-18\\
7.7	-6.17702566359865e-18\\
7.8	-3.71954108258998e-18\\
7.9	-1.46444803931838e-19\\
8	8.14129849350004e-19\\
8.1	4.64383846762686e-18\\
8.2	7.18494819196605e-18\\
8.3	8.61336845265647e-18\\
8.4	9.92999785152978e-18\\
8.5	1.58416838195422e-17\\
8.6	2.66983733060913e-17\\
8.7	4.93942768621811e-17\\
8.8	4.22786430711668e-17\\
8.9	-3.05548923324089e-17\\
9	-7.08633319501743e-16\\
9.1	-2.69647760263627e-15\\
9.2	2.60851269428862e-05\\
9.3	1.5072978640443\\
9.4	2.34689444500055\\
9.5	2.73038828600537\\
9.6	2.83792435113787\\
9.7	2.81045807101945\\
9.8	2.74672445728817\\
9.9	2.70237395243616\\
10	2.68929797954738\\
};
\addlegendentry{$u_2$};
\end{axis}
\end{tikzpicture}%
        \caption{$\eps=10^{-4}$}
    \end{subfigure}

    \bigskip 

    \begin{subfigure}[t]{0.495\textwidth}
        \begin{tikzpicture}

\begin{axis}[%
width=\textwidth,
xmin=0,
xmax=10,
ymin=-5,
ymax=30,
legend style={legend pos=north west,legend style={draw=none}},
]
\addplot [color=mycolor1,solid,line width=1.5pt]
  table[row sep=crcr]{%
0	-1.22153548928988e-17\\
0.1	-2.58314522831958e-18\\
0.2	1.45699941254591e-18\\
0.3	5.88288444500665e-18\\
0.4	1.83761131686626e-17\\
0.5	3.24571513956003e-17\\
0.6	8.21842012826903e-18\\
0.7	-3.17135065389612e-17\\
0.8	-9.14246094027283e-19\\
0.9	1.19023710267352e-17\\
1	-1.04349074990762e-17\\
1.1	1.27139456871494e-17\\
1.2	-4.93031725279912e-17\\
1.3	1.38831095061028e-05\\
1.4	8.40435852488407\\
1.5	12.788429186744\\
1.6	14.831898793176\\
1.7	15.8037278251043\\
1.8	16.4607370229118\\
1.9	17.14405303751\\
2	17.92866413523\\
2.1	18.75463660339\\
2.2	19.5166863622522\\
2.3	20.1139239385064\\
2.4	20.4710068528401\\
2.5	20.5426074777727\\
2.6	20.3101809342333\\
2.7	19.7764579767769\\
2.8	18.9602430944558\\
2.9	17.892305225973\\
3	16.6122573521621\\
3.1	15.1660342720414\\
3.2	13.6036086424512\\
3.3	11.9767409962055\\
3.4	10.3367191341895\\
3.5	8.73216163480294\\
3.6	7.20701555735674\\
3.7	5.79887860274081\\
3.8	4.53771495154512\\
3.9	3.44492882794764\\
4	2.53260948636277\\
4.1	1.80261829056523\\
4.2	1.24510461423067\\
4.3	0.836182588367379\\
4.4	0.535023205253864\\
4.5	0.281846371723343\\
4.6	-4.79851731068616e-15\\
4.7	-4.61182337241283e-15\\
4.8	-2.53693353968124e-15\\
4.9	-1.43089112301449e-15\\
5	-1.03838177086122e-15\\
5.1	-6.84907222365276e-16\\
5.2	-4.18492501564775e-16\\
5.3	-2.4162323919986e-16\\
5.4	-1.36382765284246e-16\\
5.5	-7.69157703071632e-17\\
5.6	-4.78634151362687e-17\\
5.7	-3.36404615204775e-17\\
5.8	-2.88203061567252e-17\\
5.9	-6.8287220150072e-17\\
6	6.6951700246181\\
6.1	11.0829955072144\\
6.2	12.6524832084689\\
6.3	9.84534183635299\\
6.4	-1.0054849936373e-17\\
6.5	-2.04987944625448e-19\\
6.6	1.29275109988861e-18\\
6.7	1.49983373244829e-18\\
6.8	4.77562363338294e-18\\
6.9	17.3660477928895\\
7	25.7985052863818\\
7.1	28.8288670739199\\
7.2	29.1757119661304\\
7.3	28.514157315815\\
7.4	27.665358361947\\
7.5	26.9086803330632\\
7.6	26.2587903810269\\
7.7	25.6528144035108\\
7.8	25.042912452425\\
7.9	24.4075640040796\\
8	23.6951600342697\\
8.1	22.7113926508751\\
8.2	20.9702358621163\\
8.3	17.5549185929589\\
8.4	11.0908753802673\\
8.5	-5.01129972954415e-17\\
8.6	7.7986021402401\\
8.7	10.086799891966\\
8.8	7.41242662236548\\
8.9	-1.09851237613068e-16\\
9	-5.59074973245512e-17\\
9.1	-7.97182454863935e-17\\
9.2	-1.49271888053408e-16\\
9.3	-3.54504112936121e-16\\
9.4	-1.90098489185661e-15\\
9.5	1.12447668751981\\
9.6	1.31134697367383\\
9.7	1.10707594971328\\
9.8	0.851533622115685\\
9.9	0.699701031230265\\
10	0.658802448488668\\
};
\addlegendentry{$u_1$};
\addplot [color=mycolor2,solid,line width=1.5pt]
  table[row sep=crcr]{%
0	10.5937164567864\\
0.1	10.5043038970764\\
0.2	10.2855604487125\\
0.3	10.018992027441\\
0.4	9.76466557857828\\
0.5	9.58473790550999\\
0.6	9.55484922833545\\
0.7	9.75310690234282\\
0.8	10.2145898823418\\
0.9	10.8390522681001\\
1	11.2477180160489\\
1.1	10.60966188439\\
1.2	7.50847588613408\\
1.3	1.39789766169244e-05\\
1.4	-3.19642881828911e-17\\
1.5	-3.28177770946269e-18\\
1.6	-1.74279531795022e-17\\
1.7	-4.19743380825601e-18\\
1.8	-1.05831402188195e-17\\
1.9	-1.98928006228261e-17\\
2	-7.08871067226482e-18\\
2.1	1.3006520529491e-18\\
2.2	1.01575356751482e-18\\
2.3	-7.26909124547488e-19\\
2.4	-1.22543468590029e-18\\
2.5	-7.92146066591385e-19\\
2.6	-6.20777786675856e-19\\
2.7	9.70462138872109e-20\\
2.8	1.27292707533486e-18\\
2.9	2.67319637044301e-18\\
3	5.41744276553932e-18\\
3.1	8.8221976171214e-18\\
3.2	1.33188552090633e-17\\
3.3	1.86401197186999e-17\\
3.4	2.43260447932451e-17\\
3.5	3.10987976546726e-17\\
3.6	3.81928392600952e-17\\
3.7	4.59938979110481e-17\\
3.8	5.69859520058178e-17\\
3.9	7.293035180553e-17\\
4	6.79700361356867e-17\\
4.1	1.71716544621773e-17\\
4.2	-1.93754152518405e-16\\
4.3	-1.39575785643225e-15\\
4.4	-4.44039707391974e-15\\
4.5	-2.45207241903746e-14\\
4.6	0.373158024057704\\
4.7	0.620725877015849\\
4.8	0.808378356321484\\
4.9	0.991009971353424\\
5	1.21399420203202\\
5.1	1.5230599302843\\
5.2	1.97494502438037\\
5.3	2.64062301863278\\
5.4	3.59097358861914\\
5.5	4.85236471752595\\
5.6	6.31985124821012\\
5.7	7.6242332113239\\
5.8	7.97305616980135\\
5.9	6.03406314179284\\
6	-7.17798959746587e-17\\
6.1	-1.3323505083438e-17\\
6.2	-1.17948794493325e-17\\
6.3	-1.18801949610268e-17\\
6.4	12.744338854985\\
6.5	20.5939128508775\\
6.6	24.7670551886493\\
6.7	24.8504767743967\\
6.8	18.1922290306783\\
6.9	8.46463141929813e-18\\
7	2.41056424254746e-18\\
7.1	3.48416537437199e-18\\
7.2	3.76393730214379e-18\\
7.3	4.78459795711217e-18\\
7.4	5.61389872221631e-18\\
7.5	6.03515213468727e-18\\
7.6	5.83097023502943e-18\\
7.7	5.23121061571231e-18\\
7.8	4.33715575966547e-18\\
7.9	3.07930670674333e-18\\
8	2.62076755699923e-18\\
8.1	2.28127464484642e-18\\
8.2	2.54193520313867e-19\\
8.3	6.97673122819703e-19\\
8.4	-7.78618549378509e-18\\
8.5	8.65057367578293\\
8.6	-7.96877967354631e-17\\
8.7	-3.58264642161549e-17\\
8.8	-1.08751398861147e-16\\
8.9	7.00376864724855\\
9	8.60187890248613\\
9.1	7.51790167095969\\
9.2	5.51544996821358\\
9.3	3.47739173660612\\
9.4	1.66495304312778\\
9.5	-3.04958811032094e-15\\
9.6	-1.77983958492714e-15\\
9.7	-2.00254386318283e-15\\
9.8	-2.43358884876628e-15\\
9.9	-2.01809939941813e-15\\
10	-2.42975359495372e-15\\
};
\addlegendentry{$u_2$};
\end{axis}
\end{tikzpicture}%
        \caption{$\eps=10^{-5}$}
    \end{subfigure}
    \hfill
    \begin{subfigure}[t]{0.495\textwidth}
        \begin{tikzpicture}

\begin{axis}[%
width=\textwidth,
xmin=0,
xmax=10,
ymin=-10,
ymax=60,
legend style={legend pos=north west,legend style={draw=none}},
]
\addplot [color=mycolor1,solid,line width=1.5pt]
  table[row sep=crcr]{%
0	2.16149960883215e-15\\
0.1	1.28591215466659e-15\\
0.2	1.39040894337162e-15\\
0.3	1.13795470684583e-15\\
0.4	8.5951954545362e-16\\
0.5	3.52905047192241e-16\\
0.6	-1.49546653128139e-16\\
0.7	-1.20737578869597e-15\\
0.8	-2.68630441095681e-15\\
0.9	-4.32647437615509e-15\\
1	-4.23343387295485e-15\\
1.1	-2.97128838748181e-15\\
1.2	-2.74509105920281e-15\\
1.3	10.0929487334449\\
1.4	11.3651038101749\\
1.5	11.1300319522726\\
1.6	11.7228583984751\\
1.7	13.0853969671263\\
1.8	14.6956501223422\\
1.9	16.2304494294554\\
2	17.574842641927\\
2.1	18.7000248942193\\
2.2	19.5892111905222\\
2.3	20.2199005938811\\
2.4	20.5688163358224\\
2.5	20.6196022848503\\
2.6	20.3667963148684\\
2.7	19.8169823769584\\
2.8	18.9883627093302\\
2.9	17.9097376306197\\
3	16.6189884521222\\
3.1	15.1611888833357\\
3.2	13.5863386499773\\
3.3	11.9468663765212\\
3.4	10.2950303946671\\
3.5	8.68035753140006\\
3.6	7.14727952921728\\
3.7	5.73308845015846\\
3.8	4.46655221613895\\
3.9	3.36753351508644\\
4	2.44801555334364\\
4.1	1.71323444852536\\
4.2	1.15868661955114\\
4.3	0.754571797867313\\
4.4	0.417558048658725\\
4.5	1.00372829795171e-12\\
4.6	0.377402808136494\\
4.7	0.512962827529314\\
4.8	6.28388550007546e-13\\
4.9	3.35651825082017e-13\\
5	1.84704929239799e-13\\
5.1	1.06881694433947e-13\\
5.2	3.72156778637169e-14\\
5.3	1.28652498826768e-14\\
5.4	-2.93152271794224e-15\\
5.5	-1.57696469010002e-15\\
5.6	-1.67866699802013e-15\\
5.7	7.32325853856327e-16\\
5.8	1.31714438398489e-15\\
5.9	2.02931730076507e-15\\
6	4.31693443960664e-06\\
6.1	8.29550441166764\\
6.2	4.92460911119082e-16\\
6.3	3.33020776710058e-16\\
6.4	13.7553591376071\\
6.5	21.058780309001\\
6.6	21.2656164820718\\
6.7	2.72564039004907e-05\\
6.8	23.3036307680246\\
6.9	3.07431387285885e-05\\
7	27.9444991813054\\
7.1	32.1611358517578\\
7.2	31.5889119931365\\
7.3	32.1002073818712\\
7.4	27.043632744556\\
7.5	2.8462901558146e-05\\
7.6	31.1785062088454\\
7.7	41.1252445836647\\
7.8	37.2885996862339\\
7.9	-1.65170553617111e-16\\
8	1.90230819129229e-05\\
8.1	1.40083494697912e-05\\
8.2	31.9695306062694\\
8.3	1.27341183938875e-05\\
8.4	22.6870122224918\\
8.5	19.2978038573063\\
8.6	1.05556107168929e-16\\
8.7	5.74334578626413e-16\\
8.8	7.40848309732223\\
8.9	6.98690377910946\\
9	4.33318693293874\\
9.1	1.10966842417919e-14\\
9.2	-1.68098093071053e-15\\
9.3	-1.83389852487716e-14\\
9.4	1.11356892048188\\
9.5	0.788310983851612\\
9.6	-5.98066333635379e-13\\
9.7	-7.23551022411956e-13\\
9.8	-1.49560825109661e-12\\
9.9	-2.31983198522379e-12\\
10	-3.59708032271088e-12\\
};
\addlegendentry{$u_1$};
\addplot [color=mycolor2,solid,line width=1.5pt]
  table[row sep=crcr]{%
0	10.5826539723021\\
0.1	10.4272640773983\\
0.2	10.1378101685783\\
0.3	9.89674039056706\\
0.4	9.67632890504408\\
0.5	9.41265421876403\\
0.6	9.04794564053952\\
0.7	8.55099147995029\\
0.8	8.02745124490652\\
0.9	7.94058152072529\\
1	9.16989695427504\\
1.1	12.0289423859681\\
1.2	13.0289991872873\\
1.3	-4.97820892801742e-15\\
1.4	-3.29612689778918e-15\\
1.5	-2.0538946798978e-15\\
1.6	-8.12365390202723e-15\\
1.7	4.16186985587007e-15\\
1.8	-2.6253245530543e-15\\
1.9	1.22342237464338e-15\\
2	-6.70910621751213e-16\\
2.1	-1.41627723685186e-16\\
2.2	7.73911356861214e-16\\
2.3	9.91211887867497e-16\\
2.4	6.78401489319129e-16\\
2.5	1.33552567159838e-16\\
2.6	2.52751298826972e-16\\
2.7	6.13363611122803e-16\\
2.8	7.07411898189763e-16\\
2.9	4.45471742865924e-16\\
3	-5.78106092054229e-17\\
3.1	-7.52101205446853e-16\\
3.2	-8.05019016245805e-16\\
3.3	-5.93101087344132e-16\\
3.4	-1.43352869352092e-15\\
3.5	-3.28840246694574e-15\\
3.6	-5.22226343003942e-15\\
3.7	-6.61576085547328e-15\\
3.8	-5.56155520526017e-15\\
3.9	4.80853912354255e-15\\
4	-4.54076432950378e-15\\
4.1	2.30922475143505e-14\\
4.2	1.48600000753435e-13\\
4.3	9.15056839282555e-13\\
4.4	3.45540043623713e-12\\
4.5	0.781637029665365\\
4.6	5.04358566467284e-12\\
4.7	1.89669839232535e-12\\
4.8	0.943381361152138\\
4.9	1.20648504839479\\
5	1.3630114529508\\
5.1	1.61681565876461\\
5.2	1.96814179992102\\
5.3	2.36325391574991\\
5.4	2.7732598495721\\
5.5	3.27035919494663\\
5.6	4.12317600961065\\
5.7	5.73812331789422\\
5.8	7.97053991960986\\
5.9	8.32435198465514\\
6	4.82715930576711e-06\\
6.1	2.51243410629302e-15\\
6.2	25.4981400921523\\
6.3	34.1804506798724\\
6.4	8.83212295912288e-16\\
6.5	1.32269810961016e-15\\
6.6	-3.96010839430012e-16\\
6.7	8.97741799417907e-06\\
6.8	-1.70509944383698e-16\\
6.9	9.36431019269308e-06\\
7	-1.54240285976295e-16\\
7.1	2.66105360704042e-16\\
7.2	-2.24858390739064e-16\\
7.3	1.83743447622763e-17\\
7.4	-3.96109789680189e-16\\
7.5	9.31971063504454e-06\\
7.6	-6.2224545149763e-16\\
7.7	-1.8807980338874e-16\\
7.8	-2.44734429488921e-16\\
7.9	54.7811890053383\\
8	1.06997630234809e-05\\
8.1	1.38450197498712e-05\\
8.2	-3.63574520320533e-16\\
8.3	1.25258893572933e-05\\
8.4	-3.8429035576036e-16\\
8.5	1.44775882660799e-16\\
8.6	17.8216723461223\\
8.7	13.6120917631028\\
8.8	7.0366344961641e-15\\
8.9	8.68881007500251e-15\\
9	2.23607848132936e-14\\
9.1	3.71437821095184\\
9.2	3.35794885138045\\
9.3	2.01631081973882\\
9.4	-7.35540849430984e-14\\
9.5	-3.36699445581706e-13\\
9.6	0.716329650550036\\
9.7	0.614074972564706\\
9.8	0.356962623979836\\
9.9	0.202902258675609\\
10	0.162207871021228\\
};
\addlegendentry{$u_2$};
\end{axis}
\end{tikzpicture}%
        \caption{$\eps=10^{-6}$}
    \end{subfigure}
    \caption{Optimal switching controls for different values of $\eps$}
    \label{fig:oc_diff_eta}
\end{figure}

A quantitative illustration of the dependence on $\eps$ is given in \cref{tab:results_diff_eta}, which shows that the number $N_{\mathrm{sw}}$ decreases monotonically to zero as $\eps$ increases.
At the same time, the tracking is obviously increased, resulting in a larger optimal functional value $\bar{J}$.
The fourth column confirms the switching property of the optimal control for all $\eps$.
In particular, the last row demonstrates that the proposed approach is feasible for computing switching controls for very small control cost and regularization parameters ($\alpha=10^{-6}$, $\eps=10^{-7}$, $\gamma=10^{-9}$) compared to the switching penalty ($\beta=10^{3}$).
Again, in all cases $\beta_{\max}> \alpha$, and hence the problem is genuinely nonconvex.
\begin{table}[t]
    \caption{Results for different values of $\eps$ (Example 2): optimal functional value $\bar{J}$, number of switching points $N_{\mathrm{sw}}$, switching error $\swerror$, residual norm of the optimality system \eqref{eq:opt_gamma}, final switching penalty $\beta_{\max}$}\label{tab:results_diff_eta}
    \centering
    \begin{tabular}{llclll}
        \toprule
        $\eps$ & $\bar{J}$ &$N_{\mathrm{sw}}$ & $\swerror$ & optimality & $\beta_{\max}$  \\
        \midrule
        $10^{-2}$ & $2.088$ & $0$ & $5\cdot 10^{-17}$  & $2\cdot 10^{-11}$ & $10^{-2}$ \\ 
        $10^{-3}$ & $1.202$ & $1$ & $4\cdot 10^{-22}$  & $1\cdot 10^{-11}$ & $10^0 $ \\ 
        $10^{-4}$ & $0.542$ & $4$ & $6\cdot 10^{-8 }$  & $6\cdot 10^{-11}$ & $10^3 $ \\ 
        $10^{-5}$ & $0.318$ & $9$ & $2\cdot 10^{-10}$  & $2\cdot 10^{-11}$ & $10^3 $ \\
        $10^{-6}$ & $0.274$ & $14$& $3\cdot 10^{-10}$  & $2\cdot 10^{-12}$ & $10^3 $ \\  
        $10^{-7}$ & $0.124$ & $22$& $2\cdot 10^{-11}$  & $1\cdot 10^{-12}$ & $10^3 $ \\  
        \bottomrule
    \end{tabular}
\end{table}

\section{Conclusion}

Penalization of switching constraints leads to a nonconvex optimal control problem if the switching penalty is larger than the control costs. Under additional $H^1$ regularization or restriction to finite-dimensional controls, existence of solutions can be shown. Using tools from nonsmooth analysis allows deriving optimality conditions and showing (for finite-dimensional controls) exact penalization properties. These optimality conditions are amenable to numerical solution via a (still nonconvex) Moreau--Yosida regularization and a semismooth Newton method.

By virtue of the embedding of $H^1(0,T)$ into $C([0,T])$, switching of optimal controls can only occur in points where $\bar u_1(t) = \bar u_2(t)=0$.
This is not the case if $H^1$ regularization is replaced with regularization in the space of functions of bounded variation; however, this would introduce new difficulties due to the additional nonsmoothness and the more complicated functional-analytic setting and is therefore the subject of further work.

It is possible to extend the presented approach to consider additional control constraints. Necessary optimality conditions similar to \eqref{eq:opt_incl} involving the classical normal cone of the (convex) constraint set can be obtained from \cite[Proposition 10.36]{Clarke:2013}. For the case of simple coordinate-wise box constraints, one can then apply the now classical semismooth Newton method as in, e.g., \cite{HIK:2002}. For the corresponding treatment of polygonal constraints on vector-valued controls, one can partially rely on results from \cite{KunischLu:2013} to develop semismooth Netwon methods.
Finally, another interesting aspect for future consideration would be to allow for more general switching constraints of the form $u(t) \in U$, where $U$ is a cone generated by a finite set in $L^2(0,T;\R^N)$.

\section*{Acknowledgments}

This work was supported in part by the Austrian Science Fund (FWF) under grant SFB {F}32 (SFB ``Mathematical Optimization and Applications in Biomedical Sciences''), by BioTechMed-Graz, and by the ERC advanced grant 668998 (OCLOC) under the EU's H2020 research program.

\printbibliography

\end{document}